\documentclass[preprint,12pt,nopreprintline]{elsarticle}

\usepackage{amsmath,amssymb,amsthm,mathtools,mathrsfs,bm,microtype}
\usepackage{hyperref}
\setcitestyle{numbers,square,comma,sort&compress}
\hypersetup{colorlinks=true,linkcolor=blue,citecolor=blue,urlcolor=blue}
\allowdisplaybreaks[3]
\numberwithin{equation}{section}
\newtheorem{theorem}{Theorem}[section]
\newtheorem{lemma}[theorem]{Lemma}
\newtheorem{proposition}[theorem]{Proposition}
\newtheorem{corollary}[theorem]{Corollary}
\theoremstyle{definition}
\newtheorem{definition}[theorem]{Definition}
\theoremstyle{remark}
\newtheorem{remark}[theorem]{Remark}
\newcommand{\vct}[1]{\boldsymbol{#1}}
\newcommand{\mtrx}[1]{\mathbf{#1}}
\newcommand{\Id}{\mtrx{I}}
\newcommand{\R}{\mathbb R}
\newcommand{\Om}{\Omega}
\newcommand{\NN}{\mathcal N}
\newcommand{\BB}{\mathcal B}
\newcommand{\QQ}{\mathcal Q}
\newcommand{\A}{\mathcal A}
\newcommand{\Pn}{\mtrx{P}}
\newcommand{\Qn}{\mtrx{Q}}
\DeclareMathOperator{\Dom}{D}
\DeclareMathOperator{\spec}{\sigma}
\DeclareMathOperator{\divv}{div}
\usepackage{fancyhdr}
\newcommand{\BurgersRunningTitle}{Spectral stability of negative-far-field Burgers waves}
\fancypagestyle{burgersrunning}{%
  \fancyhf{}%
  \fancyhead[L]{\small\ifodd\value{page}\relax\itshape\BurgersRunningTitle\else\thepage\fi}%
  \fancyhead[R]{\small\ifodd\value{page}\relax\thepage\else\itshape Q. Zhao\fi}%
}
\begin{document}
\begin{frontmatter}
\title{Spectral stability and nonlinear asymptotics of radial steady states for the viscous vector Burgers equation in an exterior domain with negative far field}
\author{Qingsong Zhao}
\affiliation{organization={School of Mathematics and Science, Nanyang Institute of Technology}, city={Nanyang},postcode={473004},country={China}}
\ead{qqsszhao@nyist.edu.cn}
\begin{abstract}
We study radial stationary waves of the multidimensional viscous vector Burgers equation in the exterior of a ball, with strictly negative radial far-field data. For each fixed wave satisfying a strict modified-Bessel admissibility criterion, a positive-weight conjugation yields a self-adjoint linearized operator with a positive isolated lowest eigenvalue below its essential spectrum. We establish global nonlinear stability under general continuous vector perturbations that are small in a mixed weighted-$L^2$ and unweighted-$L^\infty$ norm, without requiring small initial gradients or compact support. The lowest-eigenvalue exponential decay rate is attained for a nonempty class of admissible data. We further identify the nonlinear lowest-mode amplitude and derive refined weighted-$L^2$ and physical-variable $L^\infty$ remainder estimates, including the possible resonant case, together with a quadratic expansion in the initial data with a uniform cubic remainder.
\end{abstract}
\end{frontmatter}

\section{Introduction and main results}\label{sec:introduction}

\subsection{Background and related work}
Let $n\ge2$, $\mu>0$, and $r_0>0$, and denote the exterior of the ball
$B_{r_0}$ by $\Om=\R^n\setminus\overline{B_{r_0}}$. Set $r=|\vct{x}|$ and $\vct{\omega}=\vct{e}_r=\vct{x}/r\in\mathbb S^{n-1}$. We study the viscous vector Burgers
system
\[
 \vct{u}_t+(\vct{u}\cdot\nabla)\vct{u}=\mu\Delta\vct{u},
 \qquad (t,\vct{x})\in(0,\infty)\times\Om,
\]
with radial boundary data and a strictly negative radial far field,
$\vct{u}(t,r_0\vct{\omega})=v_-\vct{\omega}$ and
$\vct{u}(t,\vct{x})\to-b\vct{e}_r$ as $r\to\infty$, where $b>0$ and $v_+=-b<0$.
Norms and the functional-space conventions used below are specified in the Notation subsection.

This nonlinear parabolic system describes the interaction of convection
and diffusion; its classical formulation goes back to Burgers~\cite{Burgers1948}. In one space dimension,
the transformations of Hopf~\cite{Hopf1950} and Cole~\cite{Cole1951} reduce
the scalar equation to the heat equation. Higher-dimensional Cole--Hopf-type
transformations are available for particular potential formulations
~\cite{Vaganan2012}, but they do not apply to general vector perturbations with
arbitrary angular dependence. The present problem therefore requires a genuinely
vector-valued evolution analysis.

The exterior-domain theory for multidimensional Burgers equations has been
developed primarily for radially symmetric solutions. Hashimoto
~\cite{Hashimoto2014,Hashimoto2016} studied radial large-time behavior and
boundary regimes associated with rarefaction waves. Hashimoto and Matsumura
~\cite{HashimotoMatsumura2019} analyzed convergence toward nonlinear waves,
including stationary and composite profiles. Yang, Zhao, and Zhao
~\cite{YangZhaoZhao2021} established radial stability and decay results in
higher-dimensional exterior domains, while Zhao and Zhao~\cite{ZhaoZhao2021}
proved a necessary-and-sufficient existence criterion and radial stability
result in two dimensions. These works identify the stationary-wave branches
and their radial dynamics, but do not resolve the spectral evolution of general
nonradial vector perturbations.

For general multidimensional perturbations, Hashimoto~\cite[Theorem~2.3]{Hashimoto2020}
proved asymptotic stability in exterior domains for $n\ge3$ under the zero
far-field condition $v_+=0$ and the boundary restriction
$v_-<2\mu(n-2)/r_0$. That result uses smallness of the stationary flow in
$L^n$, of its gradient in $L^{n/2}$, and of the initial perturbation in a
$W^{1,n/2}$-type norm. The argument is related to stability methods for
stationary Navier--Stokes flows in exterior domains; see Borchers and
Miyakawa~\cite{BorchersMiyakawa1990} and Kozono and Ogawa~\cite{KozonoOgawa1994}.
The present paper treats the strictly negative far field $v_+=-b<0$ in every
dimension $n\ge2$, for each fixed stationary wave in the strict Bessel-admissible
range. Its initial-data condition is instead the mixed smallness
$\|v_0\|_\infty+\|hv_0\|_2\ll1$, and no initial-gradient assumption is made.
The two data classes are different; no inclusion between them is asserted.

The use of spectral information to describe nonlinear long-time behavior is
well established in the theory of evolution equations. Henry~\cite{Henry1981}
provides a semigroup framework for nonlinear stability near equilibria. Foias
and Saut~\cite{FoiasSaut1984,FoiasSaut1987} developed spectral-manifold and
normal-form methods for Navier--Stokes-type evolution, while Schonbek
~\cite{Schonbek1980} studied decay for parabolic conservation laws. Related
questions about leading modes and possible decay rates for semilinear parabolic
equations were considered by Ghisi, Gobbino, and Haraux
~\cite{GhisiGobbinoHaraux2015}; higher-order asymptotic expansions for periodic
Navier--Stokes flows were studied by Hoang and Martinez~\cite{HoangMartinez2017}.
These works provide the broader evolution-equation context, but do not give the
exterior-domain vector Burgers spectral structure or the nonlinear amplitude
estimates established here.

Related stability questions for radial viscous flows also arise in compressible
Navier--Stokes systems with inflow/outflow boundary conditions; see the preprint
by Feireisl, Gwiazda, and \v{S}wierczewska-Gwiazda (arXiv:2412.14577, 2024).
The governing equations and solution concepts in that setting are different
from those considered here.

\subsection{Negative far-field stationary states and spectral setting}
Fix $b>0$, and set $\kappa=b/(2\mu)$ and $\nu=(n-2)/2$. We consider the
initial-boundary value problem
\begin{equation}\label{eq:burgers}
\left\{\begin{aligned}
&\vct{u}_t+(\vct{u}\cdot\nabla)\vct{u}=\mu\Delta\vct{u},&&t>0,\ \vct{x}\in\Om,\\
&\vct{u}(t,r_0\vct{\omega})=v_-\vct{\omega},&&\vct{\omega}\in\mathbb S^{n-1},\\
&\vct{u}(t,\vct{x})+b\vct{e}_r\longrightarrow0,&&r\to\infty,\\
&\vct{u}(0,\vct{x})=\vct{u}_0(\vct{x}).&&
\end{aligned}\right.
\end{equation}
Denote the modified Bessel functions by $I_\nu$ and $K_\nu$. We assume the strict admissibility condition
\begin{equation}\label{eq:criterion}
v_-^{\mathrm{crit}}=b\frac{K_{\nu+1}(\kappa r_0)}{K_\nu(\kappa r_0)},\qquad v_-<v_-^{\mathrm{crit}}.
\end{equation}
Define
\begin{equation}\label{eq:weight}
\begin{gathered}
c=\frac{bI_{\nu+1}(\kappa r_0)+v_-I_\nu(\kappa r_0)}{bK_{\nu+1}(\kappa r_0)-v_-K_\nu(\kappa r_0)},\quad
h(r)=r^{-\nu}[I_\nu(\kappa r)+cK_\nu(\kappa r)],\\
\vct{w}=-2\mu\nabla\log h,\qquad \vct{F}=\nabla\log h,\\
H=-\mu\Delta_D+\mu\kappa^2\Id-2\mu\nabla^2\log h,\quad
\lambda_1=\inf\spec(H),
\end{gathered}
\end{equation}
Here $H$ is the Dirichlet realization in $L^2(\Om;\R^n)$. We will show that $h>0$, that $H$ is self-adjoint, and that $0<\lambda_1<\mu\kappa^2$. Writing $\vct{u}=\vct{w}+\vct{v}$, we use $\vct{p}:=h\vct{v}$ as the conjugated perturbation variable.

We write the initial datum as $\vct{u}_0=\vct{w}+\vct{v}_0$, so that
$\vct{v}_0$ is the perturbation of the stationary wave. The corresponding
perturbation equation and its natural mixed initial-data space are specified
in the main results below.
Once the spectral assertions in Lemma~\ref{lem:spectrum} are established, we use throughout the notation
\begin{equation}\label{eq:projection-notation}
\begin{aligned}
 E_1&:=\ker(H-\lambda_1\Id),&\Pn&:=\mathbf{1}_{\{\lambda_1\}}(H),\\
 \Qn&:=\Id-\Pn,&\Lambda&:=\inf\spec(H|_{E_1^\perp}).
\end{aligned}
\end{equation}
In particular, $0<\lambda_1<\Lambda\le\mu\kappa^2$ and $E_1$ is finite dimensional.

The strict inequality in \eqref{eq:criterion} selects the negative-far-field Bessel branch. For comparison with the two-dimensional criterion of Zhao and Zhao~\cite{ZhaoZhao2021}, when $n=2$ the radial Riccati substitution $\psi=\phi-\mu/r$ gives $V_-:=\psi(r_0)=v_- -\mu/r_0$. The equivalent boundary threshold is
\begin{equation}\label{eq:zhao-threshold}
a_{\mathrm{Bes}}:=b\frac{K_1\!\left(\dfrac{br_0}{2\mu}\right)}{K_0\!\left(\dfrac{br_0}{2\mu}\right)}-\frac{\mu}{r_0},\qquad
V_-<a_{\mathrm{Bes}}\ \Longleftrightarrow\ v_-<v_-^{\mathrm{crit}}\quad(n=2).
\end{equation}
Zhao and Zhao~\cite[Theorem~1.1]{ZhaoZhao2021} proved a necessary-and-sufficient existence criterion $V_-<a_*$ for the same transformed boundary value and the same negative far field $v_+=-b$. Proposition~\ref{prop:steady} identifies their threshold explicitly as $a_*=a_{\mathrm{Bes}}$ (see Remark~\ref{rem:exactthreshold}); the two-dimensional existence classification itself was already established in that work.

For the negative-far-field stationary wave, the Bessel conjugation transforms
the linearized vector evolution into a self-adjoint Dirichlet operator
\[
 H=-\mu\Delta_D+\mu\kappa^2 I-2\mu\nabla^2\log h,
 \qquad \kappa=\frac{b}{2\mu},
\]
where $h$ is the positive modified-Bessel weight defined in~\eqref{eq:weight}. The resulting
linear evolution is governed by the semigroup $S(t)=e^{-tH}$. The operator has
essential spectrum $[\mu\kappa^2,\infty)$, a positive isolated lowest
eigenvalue $\lambda_1$, and infinitely many discrete eigenvalues below the
essential-spectrum threshold. The latter property is caused by a long-range
attractive tangential potential of order $-2\mu\kappa/r$, which dominates the
$O(r^{-2})$ kinetic cost of trial fields supported on distant spherical shells.

\subsection{Main results}
Define $C_0(\Omega;\mathbb R^n)$ to be the continuous vector fields on $\overline\Omega$ with zero boundary trace and uniform vanishing at infinity. The natural mixed initial-data space is
\begin{equation}\label{eq:data-space}
 X_h=\{v_0\in C_0(\Omega;\mathbb R^n): hv_0\in L^2(\Omega;\mathbb R^n)\},
 \qquad \|\vct v_0\|_{X_h}=\|\vct v_0\|_\infty+\|h\vct v_0\|_2.
\end{equation}
It is a Banach space: convergence in the two component norms identifies the
weighted limit on each bounded shell. Moreover,
$h_*:=\inf_\Omega h>0$ (Proposition~\ref{prop:steady}), so $X_h\subset L^2$
continuously. Neither smoothness nor compact support is required. For example,
\[
\vct v_0(x)=\delta h(r)^{-1}e^{-r}(1-e^{-(r-r_0)})e_1
\]
belongs to $X_h$ for every $\delta$ and has noncompact support when
$\delta\ne0$; its $X_h$ norm is $O(|\delta|)$.

\begin{theorem}[Global stability under vector perturbations]\label{thm:main1}
Fix $\vct w$ satisfying~\eqref{eq:criterion}. There exist
$\epsilon,C>0$ such that, for every
\[
\vct{v}_0\in X_h,\qquad D_0:=\|\vct{v}_0\|_\infty+\|h\vct{v}_0\|_2\le\epsilon
\]
the problem \eqref{eq:burgers} admits a unique global bounded mild solution.
For every $T<\infty$,
\[
\vct v=\vct u-\vct w\in C([0,T];C_0(\Om;\R^n)).
\]
On each local continuation interval $[s,s+\tau]$,
\[
\sup_{0<t\le\tau}\sqrt t\,\|\nabla\vct v(s+t)\|_\infty<\infty.
\]
The solution is classical for positive times and
\begin{equation}\label{eq:main1}
\|h(\vct{u}(t)-\vct{w})\|_2+\|\vct{u}(t)-\vct{w}\|_\infty
\le C D_0 e^{-\lambda_1t},\quad t\ge0.
\end{equation}
Here $\lambda_1=\inf\spec(H)>0$. Furthermore, $h\vct v\in C([0,T];L^2)$ for every $T<\infty$. Neither radial symmetry nor irrotationality is imposed on the perturbation.
\end{theorem}

The second theorem identifies the leading spectral term and its nonlinear amplitude. Define
\[
\NN(\vct{p})=-h^{-1}(\vct{p}\cdot\nabla)\vct{p}+h^{-1}(\vct{F}\cdot \vct{p})\vct{p},\qquad S(t)=e^{-tH}.
\]
For smooth fields $\vct{a},\vct{b}$, define the symmetric bilinear form
\begin{equation}\label{eq:bilin}
\BB(\vct{a},\vct{b})=-\frac1{2h}\left[(\vct{a}\cdot\nabla)\vct{b}+(\vct{b}\cdot\nabla)\vct{a}\right]
+\frac1{2h}\left[(\vct{F}\cdot \vct{a})\vct{b}+(\vct{F}\cdot \vct{b})\vct{a}\right].
\end{equation}
In Section~\ref{sec:asymproof} we extend its lowest-eigenspace projection to the appropriate mixed energy spaces.
For solutions of Theorem~\ref{thm:main1}, let $\vct{p}_0=h\vct{v}_0$.
\begin{definition}[Nonlinear asymptotic amplitude]\label{def:amplitude}
For a small-data solution of Theorem~\ref{thm:main1}, define
\begin{equation}\label{eq:amplitude}
\boxed{\quad
\A(\vct{p}_0)=\Pn \vct{p}_0+\int_0^\infty e^{\lambda_1s}\Pn\NN(\vct{p}(s;\vct{p}_0))\,\mathrm{d}s\in E_1.
\quad}
\end{equation}
The integral is a Bochner integral taking values in $E_1$.
\end{definition}
Set
\begin{equation}\label{eq:gt}
\mathcal G_\Lambda(t)=
\begin{cases}
e^{-\Lambda t},&\Lambda<2\lambda_1,\\
(1+t)e^{-2\lambda_1t},&\Lambda=2\lambda_1,\\
e^{-2\lambda_1t},&\Lambda>2\lambda_1.
\end{cases}
\end{equation}
The middle branch is conditional on $\Lambda=2\lambda_1$. No claim that this
equality occurs for a particular parameter choice is needed for the estimates.

\begin{theorem}[Nonlinear lowest-mode asymptotics]\label{thm:main2}
Under the assumptions of Theorem~\ref{thm:main1}, let
$\vct A_\infty=\A(\vct p_0)$ be the amplitude in
Definition~\ref{def:amplitude}, and let $\Lambda$ be as in
\eqref{eq:projection-notation}.
Then, for $t\ge0$,
\begin{equation}\label{eq:remainderL2}
\|\vct{p}(t)-e^{-\lambda_1t}\vct{A}_\infty\|_2
\le C D_0e^{-\Lambda t}+CD_0^2\mathcal G_\Lambda(t).
\end{equation}
In particular, for any
\begin{equation}\label{eq:eta}
\lambda_1<\eta<\min\{\Lambda,2\lambda_1\}
\end{equation}
we have $\|\vct{p}(t)-e^{-\lambda_1t}\vct{A}_\infty\|_2\le C_\eta D_0e^{-\eta t}$. Likewise,
\begin{equation}\label{eq:remainderLinfty}
\|\vct{v}(t)-e^{-\lambda_1t}h^{-1}\vct{A}_\infty\|_\infty
\le C D_0e^{-\Lambda t}+CD_0^2\mathcal G_\Lambda(t).
\end{equation}
For each $\eta$ satisfying \eqref{eq:eta}, we also have
\begin{equation}\label{eq:remainderH1}
\left(\int_t^{t+1}\|\vct{p}(s)-e^{-\lambda_1s}\vct{A}_\infty\|_{H^1}^2\mathrm{d}s\right)^{1/2}
\le C_\eta D_0e^{-\eta t}.
\end{equation}
Furthermore, with the bilinear form $\BB$ defined in \eqref{eq:bilin} and $S(t)=e^{-tH}$,
\begin{equation}\label{eq:main2-quadratic}
\A_2(\vct{q}_1,\vct{q}_2)=\int_0^\infty e^{\lambda_1s}\Pn
\BB(S(s)\vct{q}_1,S(s)\vct{q}_2)\,\mathrm{d}s
\end{equation}
defines a continuous symmetric bilinear map
$L^2\times L^2\to E_1$ in the sense of projected distributional pairings, and
\begin{equation}\label{eq:main2-Taylor}
\A(\vct{p}_0)=\Pn \vct{p}_0+\A_2(\vct{p}_0,\vct{p}_0)+\mathcal R_3(\vct{p}_0),
\quad \|\mathcal R_3(\vct{p}_0)\|_2\le CD_0^3.
\end{equation}
The quadratic expansion describes the amplitude as a function of the initial data; it does not assert a complete second-order expansion in time.
\end{theorem}
\begin{remark}[Topology of the amplitude expansion]
The map $\A_2:L^2\times L^2\to L^2$ is continuous.
The amplitude map $\A$ is defined only for mixed data
$\vct p_0=h\vct v_0$ satisfying
\[
\|\vct p_0\|_2+\|h^{-1}\vct p_0\|_\infty\le\epsilon,
\qquad \vct v_0\in X_h.
\]
In particular, \eqref{eq:main2-Taylor} does not assert $L^2$-Fr\'echet twice differentiability of $\A$ on an $L^2$-open neighborhood of the origin.
\end{remark}

\subsection{Notation}
We use $r=|x|$, $\omega=e_r=x/r$, and $\mathbb S^{n-1}$ for the unit sphere.
For scalar, vector, or tensor fields $f$ on $\Omega$, we use
$|f|$ for the absolute value, Euclidean norm, or Frobenius norm,
respectively. The operator norm of a matrix is denoted by $|\cdot|_{\mathrm{op}}$.
For $1\le q<\infty$ and $q=\infty$, respectively, set
\[
 \|f\|_q=\left(\int_\Omega |f(x)|^q\,dx\right)^{1/q},
 \qquad \|f\|_\infty=\operatorname*{ess\,sup}_{x\in\Omega}|f(x)|.
\]
The domain is written explicitly when a norm is taken on a different set.
We set
\[
\begin{aligned}
\langle f,g\rangle&=\int_\Omega f\cdot g\,dx,\qquad
\|f\|_{H^1}^2=\|f\|_2^2+\|\nabla f\|_2^2,\\
H^1&=H^1(\Omega;\mathbb R^n),\qquad
H_0^1=\overline{C_c^\infty(\Omega;\mathbb R^n)}^{\,H^1}.
\end{aligned}
\]
The spaces $L^p(I;X)$ are Bochner spaces, and $H^{-1}$ denotes the dual of
$H_0^1$. The notation $C_0(\Omega;\mathbb R^n)$ means continuous vector
fields on $\overline\Omega$ with zero boundary trace and uniform vanishing
at infinity. The Dirichlet operator domain is
\[
\operatorname{Dom}(H)=H^2(\Omega;\mathbb R^n)\cap H_0^1(\Omega;\mathbb R^n).
\]
We write $\sigma(H)$ for the spectrum and $\sigma_{\mathrm{ess}}(H)$ for
the essential spectrum. The symbol $I$ is the identity on $\mathbb R^n$.
The modified Bessel functions of the first and second kinds are denoted by
$I_\nu$ and $K_\nu$, respectively, with $\nu=(n-2)/2$ and
$\kappa=b/(2\mu)$. The positive weight is $h$, its logarithmic gradient is
$F=\nabla\log h$, and the conjugated perturbation is $\vct{p}=h\vct{v}$.
The eigenspace $E_1$ corresponds to $\lambda_1$; $\Pn$ and $\Qn=\Id-\Pn$ are the orthogonal spectral projections onto $E_1$ and $E_1^\perp$, respectively,
and $\Lambda$ is the bottom of the spectrum on $E_1^\perp$. The mixed
initial-data norm is $\|\vct v_0\|_{X_h}=\|\vct v_0\|_\infty+\|h\vct v_0\|_2$, and
$D_0$ denotes this quantity for the initial perturbation. The symbols
$\mathcal N$, $\mathcal B$, $\mathcal A$, and $\mathcal A_2$
refer to the nonlinear map, symmetric bilinear map, asymptotic amplitude,
and its quadratic correction defined above. The Dirichlet quadratic form
$\mathcal Q$ is specified in Section~\ref{sec:stationary-spectrum}.
Unless explicitly fixed (in particular, the Bessel coefficient $c$ in
\eqref{eq:weight}), $C,c>0$ denote generic constants that may change from
line to line. They depend on the fixed problem parameters and stationary
profile, and, where indicated, on fixed spectral or time-window parameters,
but not on time or higher-order norms of the initial perturbation.

\subsection{Main ideas, relation to the companion paper, and organization}
The first theorem gives an endpoint exponential bound in the mixed norm, while
the second identifies the asymptotic spatial profile after multiplication by
$e^{\lambda_1t}$. In particular, the remainder estimate implies
\begin{align*}
 e^{\lambda_1t}\vct{p}(t)&\longrightarrow\A(\vct{p}_0)
 &&\text{in }L^2(\Om;\R^n),\\
 e^{\lambda_1t}\vct{v}(t)&\longrightarrow h^{-1}\A(\vct{p}_0)
 &&\text{in }L^\infty(\Om;\R^n).
\end{align*}
The correction $\A_2$ captures nonlinear interactions in the initial-data
dependence of the leading amplitude; it is not a full second-order time
expansion. Since $H$ has essential spectrum, the result is not a complete
asymptotic expansion over the entire spectrum: only the isolated lowest mode
is resolved quantitatively here. We assert neither parameter-uniform constants
nor a fixed multiplicity of $E_1$ as $v_-\uparrow v_-^{\mathrm{crit}}$;
the equality $\Lambda=2\lambda_1$ remains a conditional case.

The conjugation identity is related to the classical Witten-Laplacian structure
on one-forms~\cite{Witten1982,HelfferNier2005}. The contribution is not the
abstract identity itself, but its use with the full vector Dirichlet condition
to exclude zero modes, locate the essential-spectrum threshold, and exploit
the long-range attractive tangential potential. Positive-time regularity and
the variational realization of the weighted evolution are justified using
standard parabolic and Hilbert-space methods
~\cite{LionsMagenes1972,DenkHieberPruess2003,Lieberman1996,Wood2007,FurutoIwabuchi2026}.

A companion manuscript by the author treats the positive-far-field endpoint.
It is not used in the present analysis; the present article concerns the strict
negative-far-field branch. The positive branch corresponds to the
critical boundary value $v_-^{\mathrm{crit}}$, whereas the strict inequality
$v_-<v_-^{\mathrm{crit}}$ selects the negative-far-field branch and leaves a
growing $I_\nu$ component in the Bessel representation. Although the negative
branch converges locally to the endpoint profile as the boundary value
approaches criticality, the convergence is not uniform on the whole exterior
domain because the far-field signs differ. The companion paper proves
unweighted $L^\infty$ asymptotic stability by moving-region comparison and
angular-momentum estimates. The present paper instead uses a growing positive
weight with bounded reciprocal, which allows the spectral decay of the
conjugated variable to control the physical perturbation and makes a nonlinear
lowest-mode expansion possible. The endpoint itself is excluded from the
strict negative-far-field theorems.

The proof uses three ingredients. The positive Bessel weight yields a
self-adjoint matrix Schr\"odinger operator; a vector square identity excludes a
zero mode, and tangential trial fields place its spectral bottom below the
essential spectrum. Heat-semigroup and variational arguments then provide
global stability without an initial-gradient smallness assumption. Finally,
lowest-mode projection and unit-time energy estimates yield the nonlinear
amplitude, the quadratic correction, and the remainder bounds.

Section~\ref{sec:stationary-spectrum} establishes the Bessel representation
and the discrete lowest eigenvalue. Section~\ref{sec:global-stability}
proves Theorem~\ref{thm:main1}; Section~\ref{sec:asymproof} proves
Theorem~\ref{thm:main2}. The appendix supplies the cutoff and initial-time
trace arguments needed in the weighted energy analysis.
\section{Stationary waves and the linearized spectrum}\label{sec:stationary-spectrum}
\subsection{The Bessel representation and positivity of the weight}
\begin{proposition}[Bessel criterion]\label{prop:steady}
Condition~\eqref{eq:criterion} characterizes existence of the
negative-far-field radial wave. Its unique profile is given by
\eqref{eq:weight}. Moreover,
\[
h_*>0,\qquad h^{-1}\in W^{2,\infty}(\Om),
\]
and
\[
\vct F,\vct w\in W^{1,\infty}(\Om;\mathbb R^n).
\]
\end{proposition}
\begin{proof}
For a radially symmetric stationary state $\vct{w}=\phi(r)\vct{e}_r$,
the radial equation is $\phi\phi'=\mu(\phi'+(n-1)\phi/r)'$.
Its first integral has the form
$\mu(\phi'+(n-1)\phi/r)-\phi^2/2=C_*$.
Because $\phi(r)\to-b$, the identity gives
$\phi'(r)\to(C_*+b^2/2)/\mu$. The finite limit of $\phi$ forces this
limiting derivative to be zero. Thus $C_*=-b^2/2$, and the first integral is
\begin{equation}\label{eq:riccati}
\phi'+\frac{n-1}{r}\phi=\frac{\phi^2-b^2}{2\mu},\qquad \phi(r_0)=v_-,\quad \phi(\infty)=-b.
\end{equation}
Substitution of $\phi=-2\mu h'/h$ gives $h''+(n-1)h'/r-\kappa^2h=0$, whose general solution is
$r^{-\nu}[aI_\nu(\kappa r)+dK_\nu(\kappa r)]$. The prescribed negative far field requires $a\ne0$; we normalize $a=1$.
Using
\[
\frac d{dz}(z^{-\nu}I_\nu(z))=z^{-\nu}I_{\nu+1}(z),\qquad
\frac d{dz}(z^{-\nu}K_\nu(z))=-z^{-\nu}K_{\nu+1}(z)
\]
and $I_\nu K_{\nu+1}+I_{\nu+1}K_\nu=1/z$, write $z_0=\kappa r_0$ and
$D=bK_{\nu+1}(z_0)-v_-K_\nu(z_0)$. We obtain \eqref{eq:weight}, together with
\[
I_\nu(z_0)+cK_\nu(z_0)=\frac{b}{z_0D},\qquad
\left(\frac{I_\nu}{K_\nu}\right)'(z)=\frac1{zK_\nu(z)^2}>0.
\]
If $D>0$, then $h$ is strictly positive for all $r\ge r_0$. If $D<0$, the normalized combination is negative at $r_0$ but positive at infinity, and therefore has a zero. If $D=0$, the boundary condition excludes the growing branch. Thus \eqref{eq:criterion} is exactly the strict criterion for a smooth, pole-free stationary state with negative far field. Uniqueness follows from uniqueness of the radial ODE.
The large-argument asymptotics of the modified Bessel functions give
\begin{equation}\label{eq:h-tail}
h(r)=C_he^{\kappa r}r^{-(n-1)/2}(1+O(r^{-1})),\quad
(\log h)'(r)=\kappa-\frac{n-1}{2r}+O(r^{-2}).
\end{equation}
Consequently $h_*:=\inf_\Om h>0$, the fields $\vct{F}$, $\nabla\vct{F}$, $\vct{w}$, and $\nabla\vct{w}$ are globally bounded, and $h^{-1}\in W^{2,\infty}$. Moreover, $\Delta h=\kappa^2h$ yields
$\divv\vct{F}+|\vct{F}|^2=\kappa^2$.

For the later spectral estimates, put $f=\log h$. Since
$\Delta h=\kappa^2h$, one has
\[
 f''+(f')^2+\frac{n-1}{r}f'=\kappa^2.
\]
Differentiable Bessel asymptotics give
$f''(r)=O(r^{-2})$, so $\nabla^2f$ is bounded and tends to
zero at infinity. Its radial and tangential eigenvalues are $f''(r)$
and $f'(r)/r$, respectively. Even if $c<0$, the positivity argument
above still applies when $D>0$, because $I_\nu/K_\nu$ is strictly
increasing. We use the Bessel identities and asymptotics recorded in
\cite{DLMF}.

At $D=0$, the decaying Bessel branch has far field $+b$, so this
endpoint is excluded. As $v_-\uparrow v_-^{\mathrm{crit}}$, the denominator
$D$ tends to zero from above, while the numerator of $c$ in \eqref{eq:weight}
tends to $b/(z_0 K_\nu(z_0))>0$ by the Bessel Wronskian identity. Hence
$c\to+\infty$. After division by $c$, the weight $h/c$ converges with all
radial derivatives on compact intervals to $r^{-\nu}K_\nu(\kappa r)$.
Consequently $w$ converges locally smoothly to the positive-far-field
profile, but not uniformly on $\Omega$: the supremum norm of the difference
is at least $2b$ for each strictly negative-far-field profile. For $n=3$, the identity
$K_{3/2}(z)/K_{1/2}(z)=1+1/z$ gives the consistency check
\[
 v_-^{\mathrm{crit}}=b+\frac{2\mu}{r_0}\qquad(n=3).
\]

The earlier sufficient condition $v_-\le b+\mu(n-1)/r_0$ for a negative far-field radial wave in dimensions $n\ge4$~\cite{HashimotoMatsumura2019,YangZhaoZhao2021} is a \emph{strict subrange} of the Bessel-admissible interval. To verify the comparison, put $z>0$, $\nu=(n-2)/2$, and
$g(z)=e^z\sqrt z\,K_\nu(z)>0$. The Bessel equation gives
\[
 g''-2g'=\frac{\nu^2-1/4}{z^2}g,
 \qquad
 g'(z)=-e^{2z}\int_z^\infty\frac{\nu^2-1/4}{s^2}e^{-2s}g(s)\,ds.
\]
The second identity follows by integrating from infinity and using the standard asymptotics of $K_\nu$. Thus $g'<0$ when $n\ge4$, while $g'\equiv0$ when $n=3$. The Bessel derivative identity then yields the exact comparison
\[
 \frac{K_{\nu+1}(z)}{K_\nu(z)}-\left(1+\frac{n-1}{2z}\right)
 =-\frac{g'(z)}{g(z)}.
\]
The right-hand side is strictly positive for $n\ge4$ and zero for $n=3$.
Consequently, $v_-^{\mathrm{crit}}>b+\mu(n-1)/r_0$ for $n\ge4$, and equality holds for $n=3$. This comparison concerns explicit admissibility bounds; it does not suggest that the earlier sufficient condition was presented as necessary.
\end{proof}

\begin{remark}[The two-dimensional threshold]\label{rem:exactthreshold}
Fix $n=2$ and $v_+=-b<0$. Zhao and Zhao~\cite[Theorem~1.1]{ZhaoZhao2021}
show that, with $V_-=v_--\mu/r_0$, admissibility is equivalent to
$V_-<a_*$. Proposition~\ref{prop:steady} gives the equivalent condition
$V_-<a_{\mathrm{Bes}}$ for the same radial problem. Hence the two thresholds
coincide:
\begin{equation}\label{eq:exact-old-threshold}
 a_* = b\frac{K_1(br_0/(2\mu))}{K_0(br_0/(2\mu))}-\frac{\mu}{r_0}.
\end{equation}
Thus the Bessel representation supplies an explicit formula for the previously
characterized threshold, rather than a new existence classification in dimension two.
\end{remark}

\subsection{A vector square identity and the lowest eigenvalue}
\begin{lemma}[Vector spectral structure]\label{lem:spectrum}
Let $\QQ$ be the closed quadratic form of $H$, with form domain $H_0^1(\Om;\R^n)$. Then $H$ is self-adjoint and nonnegative and satisfies the exact square decomposition \eqref{eq:squares}. Moreover,
\[
0<\lambda_1<\mu\kappa^2,\quad \spec_{\rm ess}(H)=[\mu\kappa^2,\infty),\quad
\dim E_1<\infty,\quad\Lambda>\lambda_1.
\]
We do not assume in advance that the lowest eigenvalue is simple or that its eigenfunctions are radial.
\end{lemma}
\begin{proof}
Write $\vct{v}=\vct{u}-\vct{w}$, $f=\log h$, and $\vct{p}=h\vct{v}$.
The perturbation equation and its conjugated form are
\begin{equation}\label{eq:transformed}
\begin{gathered}
\vct{v}_t-\mu\Delta\vct{v}+(\vct{w}+\vct{v})\cdot\nabla\vct{v}+(\nabla\vct{w})\vct{v}=0,\\
\vct{p}_t+H\vct{p}=\NN(\vct{p}),\quad
\NN(\vct{p})=-h^{-1}(\vct{p}\cdot\nabla)\vct{p}
+h^{-1}(\vct{F}\cdot\vct{p})\vct{p}.
\end{gathered}
\end{equation}
The matrix-valued potential in $H=-\mu\Delta_D+\mtrx{V}$ is
$\mtrx{V}=\mu\kappa^2\Id-2\mu\nabla^2f$, which is bounded and symmetric. Therefore its closed quadratic form
\[
\QQ[\vct{p}]=\mu\|\nabla\vct{p}\|_2^2+\int_\Om\langle\mtrx{V}\vct{p},\vct{p}\rangle dx
\]
is defined on $H_0^1(\Om;\R^n)$, while the operator domain is $H^2\cap H_0^1$.
Define $D_F\vct{p}=\divv\vct{p}+\vct{F}\cdot\vct{p}$ and
$C_{F,ij}\vct{p}=\partial_ip_j-\partial_jp_i-F_ip_j+F_jp_i$.
For smooth, compactly supported fields, componentwise expansion and integration by parts show that the sum of the squared ordinary divergence and curl norms equals $\|\nabla\vct{p}\|_2^2$. The terms linear in $F$ give
$-2\int\nabla^2f[\vct{p},\vct{p}]+\int(\divv\vct{F})|\vct{p}|^2$. The quadratic terms are given by
$(F\cdot p)^2+\sum_{i<j}(F_ip_j-F_jp_i)^2=|F|^2|p|^2$.
Since $\divv F+|F|^2=\kappa^2$, we obtain
\begin{equation}\label{eq:squares}
\boxed{\ \QQ[\vct{p}]/\mu=\|D_F\vct{p}\|_2^2+\sum_{i<j}\|C_{F,ij}\vct{p}\|_2^2.\ }
\end{equation}
The two squares are the formal one-form Witten differential and adjoint
norms, with $d_{-f}=e^f d e^{-f}$; see \cite{Witten1982,HelfferNier2005}.
Here the form domain is the full vector Dirichlet space $H_0^1$, not an
absolute or relative Hodge boundary realization. Boundedness of the
coefficients and density of $C_c^\infty$ extend the identity to this domain
without boundary flux terms.
In particular, $H\ge0$; moreover, if $\vct{p}=h\vct{v}$,
\begin{equation}\label{eq:square-v}
\QQ[h\vct{v}]=\mu\!\int_\Om h^2\!\left[|\divv\vct{v}+2F\cdot\vct{v}|^2+\sum_{i<j}|\partial_iv_j-\partial_jv_i|^2\right]dx.
\end{equation}
By \eqref{eq:h-tail} and the radial Hessian decomposition,
\[
\nabla^2f=f''\vct{e}_r\otimes\vct{e}_r+\frac{f'}r(\Id-\vct{e}_r\otimes\vct{e}_r)
\]
we see that $\nabla^2f\to0$ at infinity. The potential difference $\mtrx{V}-\mu\kappa^2\Id$ is a bounded multiplier vanishing uniformly at infinity. By local $H^2$--Rellich compactness and a small-norm estimate on the tail, it is relatively compact with respect to the Dirichlet Laplacian resolvent. Weyl's theorem~\cite{Kato1995} yields
\begin{equation}\label{eq:ess}
\spec_{\mathrm{ess}}(H)=[\mu\kappa^2,\infty).
\end{equation}
To exclude zero modes, suppose $H\vct{p}=0$. Then \eqref{eq:square-v} forces $\vct{v}=h^{-1}\vct{p}$ to satisfy
\[
\partial_iv_j-\partial_jv_i=0,\quad\divv(h^2\vct{v})=0,\quad \vct{v}|_{r_0}=0.
\]
Local elliptic regularity justifies the following calculations on compact spherical shells. Construct the single-valued radial primitive
$\zeta(r,\omega)=\int_{r_0}^r\vct{v}(s\omega)\cdot\omega\,ds$.
The curl-free identities and the complete zero vector trace at the boundary imply $\nabla\zeta=\vct{v}$ and $\zeta(r_0)=\partial_r\zeta(r_0)=0$; no simply connectedness assumption is needed, even for $n=2$.
Pairing $\divv(h^2\nabla\zeta)=0$ with spherical harmonics, each angular coefficient satisfies
\[
z_\ell''+\left(\frac{n-1}r+2\frac{h'}h\right)z_\ell'-\frac{\ell(\ell+n-2)}{r^2}z_\ell=0,\qquad z_\ell(r_0)=z_\ell'(r_0)=0.
\]
Uniqueness of the ODE initial-value problem and completeness of the spherical harmonics imply $\vct{p}=0$.

The exclusion of zero modes relies on the full vector Dirichlet trace, rather than merely on vanishing of the normal component $v\cdot e_r$.
To make this explicit, choose local tangential coordinates $\theta$ on the sphere. If $\partial_i v_j-\partial_j v_i=0$, then
\[
 \partial_\theta\bigl[v(s\omega)\cdot\omega\bigr]
 =\partial_s\bigl[s\,v(s\omega)\cdot\partial_\theta\omega\bigr].
\]
Integrating radially from $r_0$, the \emph{full vector} zero trace at the inner boundary ensures
$\partial_\theta\zeta(r,\omega)=r v(r\omega)\cdot\partial_\theta\omega$,
which, together with $\partial_r\zeta=v\cdot\omega$, gives $\nabla\zeta=v$.
Thus in dimension two as well, simply connectedness of the exterior domain is unnecessary, and every spherical-harmonic component has zero Cauchy data.

In the essential-spectrum argument, the product of the truncated potential multiplier and the Dirichlet Laplacian resolvent is compact on each bounded spherical shell by
$H^2$ regularity and the compact Rellich embedding. The operator norm of the exterior part is bounded by
$\sup_{r>R}|\nabla^2f|\to0$. This establishes relative compactness; we do not infer the essential-spectrum threshold merely from pointwise convergence of the matrix potential.

Combining $H\ge0$ and \eqref{eq:ess} with the exclusion of zero from the spectrum yields $\lambda_1>0$.

Next choose $A^\top=-A\ne0$ and $\chi\in C_c^\infty(1,2)$ and construct the compactly supported tangential field
$\vct{p}_R(x)=\chi(r/R)A\vct{e}_r$. This is also admissible for $n=2$; since $A\vct{e}_r\perp\vct{e}_r$, we have
\[
\frac{\QQ[\vct{p}_R]}{\|\vct{p}_R\|_2^2}
=\mu\kappa^2+\mu\frac{\|\nabla\vct{p}_R\|_2^2}{\|\vct{p}_R\|_2^2}
-2\mu\frac{\int(f'/r)|\vct{p}_R|^2}{\|\vct{p}_R\|_2^2}.
\]
For completeness, we spell out the trial-field calculation to rule out any unaccounted $R^{-1}$ term. Set
\[
 C_0=\int_1^2\chi(s)^2s^{n-1}ds>0,\quad
 C_1=\int_1^2\chi(s)^2s^{n-2}ds>0,
\]
\[
 C_2=\int_1^2\bigl(\chi'(s)^2s^{n-1}+(n-1)\chi(s)^2s^{n-3}\bigr)ds.
\]
The first spherical-harmonic mode satisfies
\[
 \int_{\mathbb S^{n-1}}|\nabla_{\mathbb S}(A\omega)|^2
 =(n-1)\int_{\mathbb S^{n-1}}|A\omega|^2,
\]
and therefore
$\|\nabla p_R\|_2^2/\|p_R\|_2^2=C_2/(C_0R^2)$.
Moreover,
\[
\frac{\int (f'/r)|p_R|^2}{\|p_R\|_2^2}
=\frac{\kappa C_1}{C_0R}+O(R^{-2}).
\]
Thus the negative leading $R^{-1}$ contribution strictly dominates the
gradient term.
Since $p_R$ is supported in $R<r<2R$, there are no boundary fluxes at either the inner boundary or infinity.

The Rayleigh--Ritz principle and \eqref{eq:ess} establish
\begin{equation}\label{eq:spectral}
\begin{aligned}
0<\lambda_1&<\mu\kappa^2,\qquad E_1:=\ker(H-\lambda_1\Id)\ne\{0\},\\
\dim E_1&<\infty,\qquad \Lambda:=\inf\spec(H|_{E_1^\perp})>\lambda_1.
\end{aligned}
\end{equation}
Furthermore, for every $p\in H_0^1$, combining $\QQ[p]\ge\lambda_1\|p\|_2^2$ with boundedness of the potential gives
\begin{equation}\label{eq:form-coercive}
\|p\|_{H^1}^2\le C_H\QQ[p].
\end{equation}

In particular, the estimates hold on the form domain: for any $p\in H_0^1(\Om;\R^n)$,
\begin{equation}\label{eq:form-coercive-2}
\|p\|_2^2\le\lambda_1^{-1}\QQ[p],\qquad
\|p\|_{H^1}^2\le C_H\QQ[p].
\end{equation}
The second inequality follows from the first and boundedness of the potential; it does not require $p\in\Dom(H)$. Consequently, the energy identities below are justified through quadratic-form pairings and the variational time chain rule, without requiring $Hp(t)$ to be defined pointwise in $L^2$.
\end{proof}

\begin{proposition}[Long-range tangential attraction]\label{prop:infinitely-many}
The operator $H$ has infinitely many eigenvalues in $(0,\mu\kappa^2)$, counted with multiplicity. Their only accumulation point is $\mu\kappa^2$.
\end{proposition}
\begin{proof}
Fix a nonzero skew-symmetric matrix $A$ and $\chi\in C_c^\infty(1,2)$ as in the proof of Lemma~\ref{lem:spectrum}. For all sufficiently large $R$, the corresponding tangential field $p_R=\chi(r/R)Ae_r$ satisfies
$\QQ[p_R]-\mu\kappa^2\|p_R\|_2^2<0$.
Choose a strictly increasing sequence $R_j$ with $R_{j+1}>3R_j$, so that the supports of $p_{R_j}$ are disjoint. Both the Dirichlet form and the $L^2$ inner product are local; therefore, for every $N$ and arbitrary coefficients $c_j$,
\[
\QQ\!\left[\sum_{j=1}^N c_jp_{R_j}\right]
 -\mu\kappa^2\left\|\sum_{j=1}^N c_jp_{R_j}\right\|_2^2
=\sum_{j=1}^N |c_j|^2\bigl(\QQ[p_{R_j}]-\mu\kappa^2\|p_{R_j}\|_2^2\bigr)<0
\]
whenever the coefficient vector is nonzero. The min--max principle gives at least $N$ eigenvalues (counted with multiplicity) strictly below the bottom of the essential spectrum. Since $N$ is arbitrary, there are infinitely many. The essential-spectrum calculation in Lemma~\ref{lem:spectrum}, together with $H\ge\lambda_1>0$, rules out any other accumulation point.
\end{proof}
\begin{remark}
This refinement describes the spectral geometry but is not required for the two main theorems. In particular, it neither implies nor rules out the parameter-dependent equality $\Lambda=2\lambda_1$.
\end{remark}

\section{Global stability and the proof of Theorem~\ref{thm:main1}}\label{sec:global-stability}
\subsection{Mild solutions and weighted energy estimates}

\begin{lemma}[Dirichlet heat estimates for bounded data]\label{lem:heat-extension}
Let
\[
S_D(t)=e^{\mu t\Delta_D}.
\]
We use the scalar spaces
\[
X=C_{\mathrm{van}}(\overline\Om),\qquad
X_0=\{a\in X:a|_{\partial\Om}=0\}.
\]
where $X$ consists of continuous functions vanishing at infinity.
The heat evolution acts componentwise on vector fields. For the exterior ball,
the following estimates hold:
\begin{equation}\label{eq:heat-extension}
\begin{gathered}
S_D(t)X\subset X_0\quad(t>0),\\
\|\nabla^jS_D(t)a\|_\infty\le C_jt^{-j/2}\|a\|_\infty,
\quad j=0,1,2,\quad 0<t\le1.
\end{gathered}
\end{equation}
For $a\in X$, $t\mapsto S_D(t)a$ is continuous from $(0,\infty)$ to $X_0$, and $t\mapsto\nabla S_D(t)a$ is continuous from $(0,\infty)$ to $X$. In general, $S_D(t)a\to a$ in the uniform norm as $t\downarrow0$ only when $a\in X_0$.
\end{lemma}
\begin{proof}
The Dirichlet heat kernel defines a positive contraction on bounded measurable data. For an exterior ball, the killed-Brownian representation and boundary regularity give, for fixed $t>0$,
\[
|S_D(t)a(x)|\le\|a\|_\infty\mathbb P_x(\tau_\Om>t),
\qquad \mathbb P_x(\tau_\Om>t)\longrightarrow0
\quad\text{as }x\longrightarrow\partial\Om.
\]
Here $\tau_\Om$ is the first exit time of Brownian motion with generator $\mu\Delta$. For each fixed positive time, the convergence is uniform near the compact smooth boundary. To control the behavior at spatial infinity, split $a$ into a bounded function supported in a large ball and a remainder with uniform norm at most $\varepsilon$. Positivity and comparison with the whole-space heat kernel bound the contribution from the compact part by $C_t\exp[-c_t(|x|-R)^2]$, while the remainder contributes at most $\varepsilon$. Interior heat-kernel continuity gives $S_D(t)a\in X_0$.

The Dirichlet semigroup acts strongly continuously on $X_0$. For smooth compactly supported data, the exterior-domain estimate of Furuto and Iwabuchi~\cite[Theorem~1.1]{FurutoIwabuchi2026} applies with $p=\infty$ and derivative orders one and two: the exterior ball has the required smooth boundary. Density in the uniform norm extends the resulting bounds to $X_0$. Since $S_D(t/2)a\in X_0$, apply these estimates to $S_D(t/2)a$ and use the contraction property on the first half-step to obtain \eqref{eq:heat-extension} for $j=1,2$. The second-derivative bound and the fact that $S_D(t)a$ vanishes at infinity imply $\nabla S_D(t)a\to0$ at infinity: on unit balls far from the boundary, apply the difference-quotient inequality
\[
|\nabla u(x)|\le Cr^{-1}\sup_{B_r(x)}|u|+Cr\|\nabla^2u\|_\infty,
\]
then first let $|x|\to\infty$ at fixed $r$ and finally let $r\downarrow0$. Positive-time boundary regularity yields continuity of the gradient up to the sphere; no vanishing gradient trace is asserted. Splitting the semigroup at a fixed positive time and using strong continuity on $X_0$ proves the two time-continuity statements. See also \cite{Ouhabaz2005} for heat-semigroup and boundary regularity facts.
\end{proof}

\begin{lemma}[Convex Kato inequality in the variational class]\label{lem:weak-kato}
Let $I=[s,T]$ and suppose that
\[
 z\in C(I;L^2(\Om;\R^n))\cap L^2(I;H_0^1(\Om;\R^n)),\qquad z_t\in L^2(I;H^{-1}(\Om;\R^n))
\]
and that $z$ solves
\[
 z_t-\mu\Delta z+B(t,x)\cdot\nabla z+A(t,x)z=f
\]
in the variational sense. Assume that $B$ and $A$ are measurable,
$\|B\|_\infty\le B_*$, $\|A\|_{\infty,{\rm op}}\le M_*$,
and $f\in L^2(I;L^2)$. For $\varepsilon>0$ put
$\beta_\varepsilon(\xi)=\sqrt{|\xi|^2+\varepsilon^2}-\varepsilon$.
Then $\beta_\varepsilon(z)\in L^2(I;H_0^1)$, and
\begin{equation}\label{eq:weak-convex-kato}
 (\partial_t-\mu\Delta+B\cdot\nabla)\beta_\varepsilon(z)
 \le M_*\bigl(\beta_\varepsilon(z)+\varepsilon\bigr)+|f|
 \quad\text{in }\mathcal D'((s,T)\times\Om).
\end{equation}
In the limit $\varepsilon\downarrow0$, the same inequality holds with
$|z|$ in place of $\beta_\varepsilon(z)$ and without the term
$M_*\varepsilon$.
\end{lemma}
\begin{proof}
The map $\beta_\varepsilon$ is convex, $\beta_\varepsilon(0)=0$,
$|D\beta_\varepsilon|\le1$, and its Hessian is bounded (for fixed
$\varepsilon>0$), continuous and positive semidefinite. The Sobolev
composition rule gives $\beta_\varepsilon(z)\in L^2(I;H_0^1)$.
For any nonnegative $\varphi\in C_c^\infty((s,T)\times\Om)$,
$D\beta_\varepsilon(z)\varphi$ is an admissible test function in
$L^2(I;H_0^1(\Om;\R^n))$.

To justify the time chain rule, mollify $z$ in time on a compact
subinterval of $(s,T)$ containing the time support of $\varphi$.
The mollifications $z_\eta$ converge strongly in $L_t^2H_x^1$ and
$\partial_t z_\eta\to z_t$ strongly in $L_t^2H_x^{-1}$ on this subinterval.
The bounded continuous Hessian and the strong gradient convergence give
$D\beta_\varepsilon(z_\eta)\varphi\to
D\beta_\varepsilon(z)\varphi$ strongly in $L_t^2H_x^1$:
after extraction of an almost-everywhere convergent subsequence,
apply dominated convergence to the product of the Hessian difference
with $\nabla z$, and use the uniform Hessian bound for the remaining
term containing $\nabla(z_\eta-z)$. Thus the classical chain rule for
$z_\eta$, followed by the limit in the dual pairing, yields
\[
 \int_s^T\langle z_t,D\beta_\varepsilon(z)\varphi\rangle\,dt
 =-\int_s^T\!\int_\Om\beta_\varepsilon(z)\varphi_t\,dx\,dt.
\]
The spatial chain rule in the variational equation now gives
\begin{align*}
 &-\int_s^T\!\int_\Om\beta_\varepsilon(z)\varphi_t
 +\mu\int_s^T\!\int_\Om\nabla\beta_\varepsilon(z)\cdot\nabla\varphi
 +\int_s^T\!\int_\Om (B\cdot\nabla\beta_\varepsilon(z))\varphi\\
 &=\int_s^T\!\int_\Om D\beta_\varepsilon(z)\cdot(f-Az)\varphi\\
 &\quad-\mu\sum_{i=1}^n\int_s^T\!\int_\Om
 D^2\beta_\varepsilon(z)[\partial_i z,\partial_i z]\varphi
 \le\int_s^T\!\int_\Om (|f|+M_*|z|)\varphi .
\end{align*}
All space--time integrals are with respect to $dx\,dt$.
Since $|z|\le\beta_\varepsilon(z)+\varepsilon$, this proves
\eqref{eq:weak-convex-kato}. In fact the stronger right-hand side
$|f|+M_*|z|$ is already integrable on the full energy space.
Finally, $0\le\beta_\varepsilon(z)\le |z|$ and the Lipschitz composition
rule yield strong $L_t^2L_x^2$ convergence to $|z|$. For the gradients,
$D\beta_\varepsilon(z)\nabla z$ converges almost everywhere to
$\nabla|z|$ (zero on $\{z=0\}$), and is bounded by $|\nabla z|$.
Dominated convergence gives strong $L_t^2H_x^1$ convergence. Passing to
the limit against $\varphi$ proves the claimed modulus inequality.
No derivative of $B$ and no artificial outer-boundary limit are needed.
\end{proof}

\begin{lemma}[Local weighted mild theory]\label{lem:localweighted}
Given $\vct v_0\in X_h$, a unique local mild solution exists.
It continues while $\|\vct v(t)\|_\infty$ stays bounded. For each $T<\infty$,
\[
\vct p=h\vct v\in C([0,T];L^2)\cap L^2(0,T;H_0^1),\qquad
\NN(\vct p)\in L^2(0,T;L^2).
\]
The variation-of-constants formula \eqref{eq:duhamel} holds in the
$L^2$-Bochner sense.
\end{lemma}
\begin{proof}
On the smooth exterior ball, the Dirichlet heat semigroup $S_D(t)=e^{\mu t\Delta_D}$ satisfies the standard short-time estimates, for $0<t\le1$,
\begin{equation}\label{eq:heat}
\|S_D(t)a\|_\infty\le\|a\|_\infty,\qquad
\|\nabla S_D(t)a\|_\infty\le C_Dt^{-1/2}\|a\|_\infty.
\end{equation}
For vector- and matrix-valued functions, interpret the spaces $X=C_{\mathrm{van}}(\overline\Om)$ and $X_0=C_0(\Om)$ componentwise. The bounded-data extension in Lemma~\ref{lem:heat-extension} applies to all
$a\in C_{\mathrm{van}}(\overline\Om)$, including inputs whose inner boundary traces are nonzero.
Thus, the short-time estimate \eqref{eq:heat} and positive-time gradient continuity do not rely on approximating nonzero-trace data by $C_c^\infty(\Om)$ in the uniform norm.

With initial time $s$, let $J=[s,s+\tau]$ and work in the Banach space
\begin{equation}\label{eq:fixedspace}
\begin{aligned}
\mathfrak X_{s,\tau}=\bigl\{v\in C(J;C_0(\Om;\R^n)):\;&
\nabla v\in C((s,s+\tau];C_{\mathrm{van}}(\overline\Om;\R^{n\times n})),\\
&\sup_{0<t\le\tau}\sqrt t\,\|\nabla v(s+t)\|_\infty<\infty\bigr\}
\end{aligned}
\end{equation}
equipped with the norm
\[
\|v\|_{\mathfrak X_{s,\tau}}=
\sup_{t\in J}\|v(t)\|_\infty+
\sup_{0<t\le\tau}\sqrt t\,\|\nabla v(s+t)\|_\infty.
\]
To verify completeness, observe that every Cauchy sequence in this norm converges in $C(J;C_0)$,
while its gradients converge uniformly in $C_tC_{\mathrm{van}}$ on each $[s+\delta,s+\tau]$.
Identification of distributional derivatives then shows that the limiting gradient belongs to the stated class.
Write $W=\|w\|_\infty$ and $M=\|\nabla w\|_{\infty,\mathrm{op}}$, and define
\[
\begin{aligned}
G(v)&=(w+v)\cdot\nabla v+(\nabla w)v,\\
\Phi(v)(s+t)&=S_D(t)v_s-\int_0^t S_D(t-r)G(v)(s+r)\,dr.
\end{aligned}
\]
For $v\in\mathfrak X_{s,\tau}$ and $r>0$, observe that
$G(v)(s+r)\in C_{\mathrm{van}}(\overline\Om)$:
indeed, $w$ and $\nabla w$ are bounded and $v,\nabla v$ vanish at spatial infinity,
but $G(v)$ generally \emph{does not} have a zero trace on the inner boundary.
Accordingly, the Duhamel term must use the extension of the semigroup from $C_{\mathrm{van}}$ to $C_0$ described above,
rather than a density argument applicable only to $C_0$ inputs.
If $\|v\|_{\mathfrak X_{s,\tau}}\le K$, then
\[
\|G(v)(s+r)\|_\infty\le K(W+K)r^{-1/2}+MK.
\]
By \eqref{eq:heat} and the identity
$\int_0^t(t-r)^{-1/2}r^{-1/2}\,dr=\pi$, we obtain
\begin{align*}
\sup_{0\le t\le\tau}\|\Phi(v)(s+t)\|_\infty
&\le\|v_s\|_\infty+2K(W+K)\sqrt\tau+MK\tau,\\
\sup_{0<t\le\tau}\sqrt t\,\|\nabla\Phi(v)(s+t)\|_\infty
&\le C_D\|v_s\|_\infty+C_D\pi K(W+K)\sqrt\tau+2C_DMK\tau.
\end{align*}
For fixed $t>0$, $\nabla S_D(t-r)G(v)(s+r)$ belongs to
$C_{\mathrm{van}}$ and admits the bound
$C(t-r)^{-1/2}(r^{-1/2}+1)$, which is integrable.
Truncate an interval of length $\varepsilon$ near each endpoint of integration, use positive-time continuity, and then
let $\varepsilon\downarrow0$. This shows that the gradient of the Duhamel term is continuous on
$(s,s+\tau]$ with values in $C_{\mathrm{van}}$.
The Duhamel term itself has $C_0$ norm $O(\sqrt t)$, so
$\Phi(v)\in\mathfrak X_{s,\tau}$ with initial trace $v_s$.

For two candidate functions, express the difference of the nonlinearities as
\[
G(v_1)-G(v_2)=
(w+v_1)\cdot\nabla(v_1-v_2)
+(v_1-v_2)\cdot\nabla v_2+(\nabla w)(v_1-v_2).
\]
The same convolution calculation yields
\[
\|\Phi(v_1)-\Phi(v_2)\|_{\mathfrak X_{s,\tau}}
\le C[(W+K)\sqrt\tau+M\tau]
\|v_1-v_2\|_{\mathfrak X_{s,\tau}}.
\]
Fix $K\ge2(1+C_D)\|v_s\|_\infty+1$ and then choose $\tau$ sufficiently small.
The map is then a strict contraction on the closed ball. The local existence time depends only on the initial supremum norm
and the fixed coefficients, not on $\|\nabla v_s\|_\infty$.

The fixed point satisfies the distributional equation. Its positive-time classical regularity follows by the localized heat-equation argument in Lemma~\ref{lem:positive-regularity} below, which uses only the local fixed-point norm and local energy uniqueness, not the weighted reconstruction or the global stability estimate.

We verify uniqueness in the full bounded mild-solution class of Theorem~\ref{thm:main1}, not merely in the fixed-point ball of \eqref{eq:fixedspace}.
Let $v_1,v_2$ be two such solutions with the same initial data and put $d=v_1-v_2$.
On any positive-time window the difference satisfies
\[
d_t-\mu\Delta d+(w+v_1)\cdot\nabla d+
[(\nabla w)+(\nabla v_2)]d=0.
\]
The bracketed term denotes the corresponding matrix acting on $d$.
For any $0<\delta<t$, vector Kato comparison with bounded drift gives
\[
\|d(t)\|_\infty\le\|d(\delta)\|_\infty
\exp\!\left[M(t-\delta)+\int_\delta^t\|\nabla v_2(r)\|_\infty\,dr\right].
\]
On the first common fixed-point window, $\|\nabla v_2(r)\|_\infty\le Cr^{-1/2}$,
so the exponential integral is bounded uniformly as $\delta\downarrow0$.
Both mild solutions attain the same initial data in $C_tC_0$, hence $d(t)=0$.
Repeating the argument on subsequent restart windows and using agreement on overlaps
yields uniqueness throughout the common interval of existence.
This step requires only boundedness and the local time-weighted gradient bound specified in the theorem for competing solutions;
it does not presume that their gradients vanish at infinity.

For the local solution, the vector Kato inequality for the modulus gives
\[
(\partial_t-\mu\Delta+(w+v)\cdot\nabla)|v|\le M|v|
\]
The scalar maximum principle on the exterior domain yields
\begin{equation}\label{eq:kato}
\|v(t)\|_\infty\le e^{M(t-s)}\|v(s)\|_\infty.
\end{equation}
This bound rules out finite-time blow-up in the supremum norm. Since the local existence time depends only on that norm,
we can restart at any $s<T_{\max}$ and thus exclude finite-time termination.
The bound on $\sqrt t\,\|\nabla v(t)\|_\infty$ on the first window follows directly from the fixed point.
On each $[t_0,T]$ with $t_0>0$, gradient continuity and uniqueness over finitely many positive-time restart windows
give $\sup_{t_0\le t\le T}\|\nabla v(t)\|_\infty<\infty$.
Consequently, the finite-time $\sqrt t$-weighted gradient bound used below holds,
independently of any weighted energy estimate that has yet to be proved.

We next establish an unweighted variational formulation. Since $h_*>0$ and $v_0\in X_h$, the initial value belongs to $L^2$. Fix $T<\infty$ and freeze $B=w+v$, $A=\nabla w$ in the equation $z_t-\mu\Delta z+B\cdot\nabla z+Az=0$, $z(0)=v_0$.
The variational solution can be constructed by Galerkin approximation on truncated spherical shells with bounded coefficients.
For fixed $T$ put $B_*=\|w+v\|_{L^\infty((0,T)\times\Omega)}$ and
$M_*=\|\nabla w\|_{L^\infty,\mathrm{op}}$. For $\psi\in H_0^1$, consider the bilinear form
\[
 a_t(z,\psi)=\mu(\nabla z,\nabla\psi)_2
 +(B\cdot\nabla z,\psi)_2+(Az,\psi)_2
\]
The form is bounded and satisfies a shifted G\aa rding inequality
\[
 a_t(z,z)+C(B_*^2/\mu+M_*+1)\|z\|_2^2
 \ge \frac\mu2\|\nabla z\|_2^2.
\]
Consequently, on each finite spherical shell the Galerkin solutions satisfy
$L_t^\infty L_x^2\cap L_t^2H_0^1$ estimates independent of the outer shell radius, and their time derivatives are bounded in
$L_t^2H^{-1}$. Weak compactness yields a variational solution on the full exterior domain;
the time chain rule ensures $z\in C_tL^2$ and strong attainment of the initial data in $L^2$.
The difference of two solutions satisfies the same energy inequality, hence uniqueness. Here the convection term is estimated directly by
Young's inequality, without integrating $B$ by parts; therefore no a priori bound on $\divv B$ is required.

To identify the mild and variational solutions without assuming that the mild solution lies in $L^2$, put $d=v-z$ and fix $a>0$, $\omega_a=e^{-a\sqrt{1+|x|^2}}$. This weight is integrable and satisfies $|\nabla\omega_a|\le a\omega_a$.
Take a radial exterior cutoff $\chi_R$, equal to $1$ for $r\le R$ and to $0$ for $r\ge2R$,
with $|\nabla\chi_R|\le C/R$. On each $[\delta,T]$ test the difference equation against $\chi_R^2\omega_a d$,
without taking a divergence of the drift, and write
\[
Y_{a,R}=\int\chi_R^2\omega_a|d|^2,\quad
G_{a,R}=\int\chi_R^2\omega_a|\nabla d|^2.
\]
Young's inequality and $|\nabla\omega_a|\le a\omega_a$ give
\[
\tfrac12Y'_{a,R}+\tfrac\mu2G_{a,R}
\le\left(\frac{(\mu a+B_*)^2}{\mu}+M_*\right)Y_{a,R}
+4\mu\!\int\omega_a|\nabla\chi_R|^2|d|^2.
\]
Since $v$ is bounded on finite intervals, $z\in C_tL^2$, and $\omega_a\in L^1$, the time integral of the last term is $O(R^{-2})$.
Apply Gronwall's inequality at fixed $R$, then pass to $R\to\infty$ to obtain
\begin{equation}\label{eq:finite-exponential-uniqueness}
 \int\omega_a|d(t)|^2\le
 e^{2[(\mu a+B_*)^2/\mu+M_*](t-\delta)}
 \int\omega_a|d(\delta)|^2.
\end{equation}
Because $z(\delta)\to v_0$ in $L^2$ and $v(\delta)\to v_0$ in $C_0$,
the right-hand side tends to zero as $\delta\downarrow0$. Thus we need neither assume $v\in L^2$ in advance nor integrate
$\|\nabla v(t)\|_\infty^2\lesssim t^{-1}$ at time zero.

Hence $v=z\in C([0,T];L^2)\cap L^2(0,T;H_0^1)$. This identification is established before testing with the unbounded weight.

The identification precedes any unbounded-weight testing; the time-weighted mild gradient bound alone would not imply an integrated gradient square at time zero.
Apply Lemma~\ref{lem:truncation-principle} with $a_N=\min\{h,N\}$. Its hypotheses hold because $v\in C_tL^2\cap L_t^2H_0^1$, the frozen drift and matrix are bounded, and $|\nabla a_N|\le\|F\|_\infty a_N$. Taking the spatial cutoff limit at fixed $N$ and then using monotone convergence yields
\begin{equation}\label{eq:weighted}
 \sup_{0\le t\le T}\|hv(t)\|_2^2+\int_0^T\|hv(t)\|_{H^1}^2dt\le C_T\|hv_0\|_2^2,
\end{equation}
where the product rule and $|\nabla h|\le\|F\|_\infty h$ control the gradient. In particular $p=hv\in L^2(0,T;H_0^1)$: the trace is zero locally and the weighted field and its gradient are square-integrable globally. The transformed equation gives
\[
 \|\NN(p)\|_2\le\|v\|_\infty(\|\nabla p\|_2+\|F\|_\infty\|p\|_2),
\]
so $\NN(p)\in L^2(0,T;L^2)$ and $p_t\in L^2(0,T;H^{-1})$. The Hilbert-triple theorem supplies strong $L^2$ continuity. Its initial trace is exactly $p_0=hv_0$, because for any compactly supported test field $\phi$, $\langle p(t),\phi\rangle=\langle v(t),h\phi\rangle\to\langle p_0,\phi\rangle$; the uniform weighted bound extends this weak convergence to all $L^2$ tests. See Remark~\ref{rem:trace-detail}. These arguments apply directly to $X_h$, without approximation by differentiable initial data or a support restriction.

For the variation-of-constants formula, first choose $p_0^{(m)}\in\Dom(H)$ and
$F^{(m)}\in C_c^\infty((0,T);L^2)$ approximating
$p_0$ and $F=\NN(p)$, respectively. The corresponding strong solutions satisfy
\[
 p^{(m)}(t)=S(t)p_0^{(m)}+\int_0^tS(t-s)F^{(m)}(s)ds.
\]
Since $\|S(t)\|_{2\to2}\le1$, we have
\[
 \sup_{0\le t\le T}\Bigl\|
\int_0^t S(t-s)(F^{(m)}-F)(s)ds\Bigr\|_2
 \le\sqrt T\,\|F^{(m)}-F\|_{L^2(0,T;L^2)}.
\]
By uniqueness of linear variational solutions, passage to $m\to\infty$ yields, for the same $p$, the
$L^2$-Bochner formula, rather than a merely formal inversion of $H$.

We have therefore obtained
\begin{equation}\label{eq:duhamel}
p(t)=e^{-tH}p_0+\int_0^t e^{-(t-s)H}\NN(p(s))\,ds,
\qquad t\ge0,
\end{equation}
The integral is a Bochner integral in $L^2$. Since $H\ge\lambda_1$,
\[
\left\|\int_0^t e^{-(t-s)H}\NN(p(s))ds\right\|_2
\le\sqrt t\,\|\NN(p)\|_{L^2(0,t;L^2)}\longrightarrow0\quad(t\downarrow0).
\]
Only $p_0\in L^2$ is used; no initial $H^1$ membership is assumed.

The constant in \eqref{eq:weighted} may depend on $T$; global decay is established separately below.
\end{proof}

\subsection{Smoothing with bounded drift and positive-time regularity}
\begin{lemma}[Short-time smoothing for nonautonomous drift]\label{lem:moser}
Let $z$ be a variational solution of the equation with bounded measurable drift $B$ and matrix-valued zeroth-order term $A$, satisfying $\|B\|_\infty\le B_*$ and $\|A\|_{\infty,{\rm op}}\le M_*$. Then for every $0<\tau\le1$,
\begin{equation}\label{eq:smoothing-statement}
\|z(s+\tau)\|_\infty\le C_{n,\mu}\tau^{-n/4}
\exp[(M_*+B_*^2/\mu)\tau]\|z(s)\|_2.
\end{equation}
The estimate holds at every fixed terminal time and does not require $\operatorname{div}B\in L^\infty$.
\end{lemma}
\begin{proof}
To make the constants independent of the restart time, consider
$z_t-\mu\Delta z+B(t,x)\cdot\nabla z+A(t,x)z=0$,
$\|B\|_\infty\le B_*$ and $\|A\|_{\infty,\mathrm{op}}\le M_*$.
Apply Lemma~\ref{lem:weak-kato}, first with $\beta_\varepsilon(z)$ and then with $\varepsilon\downarrow0$, to $q=e^{-M_*(t-s)}|z|$, so that
$ q_t-\mu\Delta q+B\cdot\nabla q\le0$ in the distributional sense.
For each $m\ge2$ test this inequality first with a bounded truncation of $m q^{m-1}$ and a compactly supported exterior cutoff. Steklov averaging justifies the time chain rule. The cutoff is removed using the $L_t^\infty L_x^2\cap L_t^2H_0^1$ estimates at the preceding exponent; the bound at $m=2$ starts the induction. Thus no \emph{a priori} finite $L^m$ norm for arbitrary $m$ is assumed. The drift is estimated directly,
$ mB_*\int q^{m-1}|\nabla q|
   =2B_*\int q^{m/2}|\nabla q^{m/2}|$,
without taking a divergence of $B$. Young's inequality gives
\begin{equation}\label{eq:moser-energy}
\frac d{dt}\|q\|_m^m+\frac{2\mu(m-1)}m\|\nabla q^{m/2}\|_2^2
\le\frac{m B_*^2}{2\mu(m-1)}\|q\|_m^m,\quad m\ge2.
\end{equation}
Here $m/[2(m-1)]\le1$ and $2(m-1)/m\ge1$.
Thus, with $\theta(t)=e^{-B_*^2(t-s)/\mu}q(t)$, we have
\begin{equation}\label{eq:theta-energy}
\frac d{dt}\|\theta\|_m^m+\mu\|\nabla\theta^{m/2}\|_2^2\le0
\quad\text{for every }m\ge2.
\end{equation}
Put $\alpha=1+2/n$, $m_k=2\alpha^k$, and
$t_k=s+\frac\tau2(1-2^{-k})$ and $I_k=(t_k,s+\tau)$.
Use time cutoffs in \eqref{eq:theta-energy}, choosing
$\eta_k=0$ for $t\le t_k$ and $\eta_k=1$ for $t\ge t_{k+1}$, with
$|\eta_k'|\le C2^{k+2}/\tau$. This gives
\[
\sup_{t\in I_{k+1}}\|\theta^{m_k/2}(t)\|_2^2+
\mu\int_{I_{k+1}}\|\nabla\theta^{m_k/2}\|_2^2dt
\le \frac{C2^{k+2}}\tau\int_{I_k}\!\int_\Om\theta^{m_k}dxdt.
\]
Extend $\theta^{m_k/2}\in H_0^1(\Om)$ by zero to $\R^n$
and apply the parabolic Sobolev--Gagliardo--Nirenberg inequality
\[
\int_I\!\int_\Om|f|^{2\alpha}\le C_n
\bigl(\sup_{t\in I}\|f(t)\|_2^2\bigr)^{2/n}
\int_I\|\nabla f(t)\|_2^2dt,
\]
In dimension $n=2$ we use the $L^4$ Gagliardo--Nirenberg inequality, since $H^1$ does not embed into $L^\infty$. It follows that
\[
\|\theta\|_{L^{m_{k+1}}(I_{k+1}\times\Om)}\le
\left(\frac{C_{n,\mu}2^{k+2}}\tau\right)^{1/m_k}
\|\theta\|_{L^{m_k}(I_k\times\Om)}.
\]
Since $\sum_{k\ge0}1/m_k=(n+2)/4$,
$\sum_{k\ge0}k/m_k<\infty$, and
$\|\theta\|_{L^2((s,s+\tau)\times\Om)}\le\sqrt\tau\|z(s)\|_2$,
the iteration limit gives $\|\theta\|_{L^\infty((s+\tau/2,s+\tau)\times\Om)}
\le C_{n,\mu}\tau^{-n/4}\|z(s)\|_2$.
First establish the bound for almost every terminal time. Strong $L^2$ continuity of the variational solution and
lower semicontinuity of the $L^\infty$ norm extend the terminal estimate to every fixed $s+\tau$.
Finally, approximate the data and bounded measurable coefficients by smooth quantities while preserving the uniform coefficient bounds,
and pass to the limit by variational energy uniqueness and strong convergence in $C_tL^2$.
Reinstating the two exponential weights gives
\begin{equation}\label{eq:smoothing}
\|z(s+\tau)\|_\infty\le C_{n,\mu}\tau^{-n/4}
\exp[(M_*+B_*^2/\mu)\tau]\|z(s)\|_2,\quad0<\tau\le1.
\end{equation}
No bound on $\nabla B$, $\divv B$, or time derivatives is used.

\end{proof}

\begin{lemma}[Classical regularity for positive times]\label{lem:positive-regularity}
Every local solution in the fixed-point class \eqref{eq:fixedspace} is classical away from its initial time, locally up to the inner boundary. This assertion uses neither an initial weighted norm nor global continuation.
\end{lemma}
\begin{proof}
We establish classical regularity at positive times independently of the stability estimate. Fix arbitrary $0<\sigma<T/3$ and
$r_0<R_1<R_2<R_3$, and set $D_j=\{r_0<r<R_j\}$.
Choose a spatial cutoff $\chi_0=1$ for $r\le R_2$, vanishing near $R_3$,
and a time cutoff $\eta_0=0$ for $t\le\sigma/4$, equal to $1$ for $t\ge\sigma/2$.
On $D_3$, the function $u_0=\eta_0\chi_0 v$ solves a heat equation with zero initial and Dirichlet data, whose right-hand side is
\[
\begin{aligned}
G_0={}&\eta_0\chi_0 G(v)+\eta_0'\chi_0 v
-2\mu\eta_0\nabla\chi_0\cdot\nabla v\\
&-\mu\eta_0(\Delta\chi_0)v,\\
G(v)={}&-(w+v)\cdot\nabla v-(\nabla w)v.
\end{aligned}
\]
By the fixed-point gradient bound, this forcing belongs to every finite $L^p$. Choose $p>n+2$.
For this fixed bounded $C^\infty$ shell, the Dirichlet Laplacian with
constant principal coefficient $\mu>0$ is parameter-elliptic and
satisfies the Dirichlet Lopatinskii--Shapiro complementing condition.
Its realization on $L^p(D_3)$ has domain
$W^{2,p}(D_3)\cap W_0^{1,p}(D_3)$ and the required maximal
$L^p$ regularity (Denk--Hieber--Pr\"uss~\cite[Part~II, Theorem~8.2]{DenkHieberPruess2003};
see also the smooth-domain specialization discussed by Wood~\cite{Wood2007}).
A positive spectral shift, if used in the abstract theorem, is removed
on the finite time interval by an exponential time weight. Consequently
the zero-initial-data problem has a strong $W_p^{2,1}$ solution.
Both this solution and $u_0$ belong to $C_tL^2\cap L_t^2H_0^1$ on the bounded shell, since the time cutoff avoids the initial time and the fixed-point gradient is bounded there;
energy uniqueness identifies them.
The parabolic Sobolev--Morrey embedding gives $C^{\alpha,\alpha/2}$ regularity of $v$ and $\nabla v$ on
$[\sigma,T]\times\{r_0\le r\le R_2\}$,
where $0<\alpha<1-(n+2)/p$.
Next choose $\chi_1=1$ for $r\le R_1$ and vanishing near $R_2$,
and $\eta_1=0$ in a short interval just after $\sigma$ and equal to $1$ for $t\ge2\sigma$.
The zero-initial-data Dirichlet heat equation for $u_1=\eta_1\chi_1v$ has a right-hand side of the same form,
now in $C^{\alpha,\alpha/2}$. The $C^\infty$ spatial boundary, constant uniformly parabolic
coefficients, zero boundary and initial data, and H\"older forcing meet
the hypotheses of the Cauchy--Dirichlet Schauder existence theorem
\cite[Chapter~V, Theorem~5.14]{Lieberman1996}. In particular, the second
time cutoff makes the forcing vanish near the artificial initial time,
so the initial--boundary compatibility condition is satisfied.
Energy uniqueness identifies that classical solution with $u_1$, yielding
$v\in C^{2+\alpha,1+\alpha/2}$ on
$[2\sigma,T]\times\{r_0\le r\le R_1\}$.

The time cutoffs give zero artificial initial data and the second forcing vanishes near the artificial initial time. Thus the corner compatibility needed for Schauder estimates is automatic. Since $\sigma$ and $R_1$ are arbitrary, this proves the assertion locally up to the true boundary; no uniformity in the outer shell radius is required.
\end{proof}

\subsection{The nonlinear energy estimate and sharp exponential decay}
\begin{proof}[Proof of Theorem~\ref{thm:main1}]
For $p\in C_tL^2\cap L_t^2H_0^1$, test \eqref{eq:transformed} with $p$;
the Hilbert-triple time chain rule applies by Lemma~\ref{lem:localweighted}.
For almost every time, $v=h^{-1}p$ is bounded and $p\in H_0^1$.
Set $g=h^{-1}|p|^2$. Then $g\in L^2$ and
\[
 |\nabla g|\le 2|v||\nabla p|+|F||v||p|\in L^2.
\]
Indeed, the chain and product rules first apply to the locally
$W^{1,1}$ expression $h^{-1}|p|^2$ and give the displayed global
$L^2$ bounds. On a bounded collar of the physical boundary $h^{-1}$
is smooth and bounded; the trace of $g$ is zero, by the zero trace
of $p$ and local boundedness of $p=hv$ (or by a bounded Lipschitz
truncation of $|p|^2$ on that collar). Hence $g\in H_0^1$.
The weak divergence identity $\int p\cdot\nabla g=-\int (\divv p)g$
is therefore valid on the full exterior domain. Using
$h^{-1}p\cdot\nabla |p|^2=p\cdot\nabla g+(F\cdot p)g$, one obtains
\[
 \langle\NN(p),p\rangle
 =\tfrac12\int_\Om(\divv p+F\cdot p)g\,dx.
\]
All terms are integrable by the above bounds. This gives the exact cubic term
\begin{equation}\label{eq:energy}
\frac12\frac d{dt}\|p(t)\|_2^2+\QQ[p(t)]
=\mathscr R[p(t)],\qquad
\mathscr R[p]=\frac12\int_\Om h^{-1}(D_Fp)|p|^2dx.
\end{equation}
The cubic term is not sign-definite. Nevertheless, \eqref{eq:squares}, the spectral lower bound, and $h^{-1}|p|^2=|v||p|$ imply
\begin{equation}\label{eq:absorb}
|\mathscr R[p]|\le\tfrac12\|D_Fp\|_2\|v\|_\infty\|p\|_2
\le\frac{\|v\|_\infty}{2\sqrt{\mu\lambda_1}}\QQ[p].
\end{equation}
Put $\delta=\sqrt{\mu\lambda_1}$. On any interval where the bootstrap bound $\|v\|_\infty\le\delta$ holds,
\begin{equation}\label{eq:half-energy}
\frac12(\|p\|_2^2)'+\frac12\QQ[p]\le0,\qquad
\|p(t)\|_2\le e^{-\lambda_1t/2}\|p_0\|_2.
\end{equation}
First fix $0<\tau\le1$ and the drift bound $\|w\|_\infty+\delta$, obtaining $C_\tau$ from \eqref{eq:smoothing}. Then choose the data sufficiently small that
$e^{M\tau}\|v_0\|_\infty<\delta/4$ and
$C_\tau h_*^{-1}\|p_0\|_2<\delta/4$.
If a first hitting time of $\delta$ existed, then for $t\le\tau$ preceding that time, \eqref{eq:kato} would give
$\|v(t)\|_\infty<\delta/4$; for $t\ge\tau$, \eqref{eq:smoothing} and \eqref{eq:half-energy} would give
\[
\|v(t)\|_\infty\le C_\tau\|v(t-\tau)\|_2
\le C_\tau h_*^{-1}e^{-\lambda_1(t-\tau)/2}\|p_0\|_2<\delta/4.
\]
This contradicts the definition of the first hitting time (continuity of the supremum norm follows from $C_tC_0$). Thus $\|v\|_\infty<\delta$ globally, and
\begin{equation}\label{eq:half-decay}
\|p(t)\|_2+\|v(t)\|_\infty\le CD_0e^{-\lambda_1t/2},\quad t\ge0.
\end{equation}

Here $\tau$ and its smoothing constant are fixed before choosing $\epsilon$, so the contradiction and its constants are uniform over the mixed small-data ball, regardless of whether the data have compact support or derivatives.

To recover the sharp exponential rate, set $a(t)=\|v(t)\|_\infty/(2\sqrt{\mu\lambda_1})$.
Equations \eqref{eq:energy}--\eqref{eq:absorb} give
$\frac12(\|p\|_2^2)'+[1-a(t)]\QQ[p]\le0$.
By \eqref{eq:half-decay},
$\int_0^\infty a(s)ds\le CD_0/\lambda_1$, and hence
\[
\|p(t)\|_2\le\|p_0\|_2
\exp\left[-\lambda_1t+\lambda_1\int_0^ta(s)ds\right]
\le CD_0e^{-\lambda_1t}.
\]
Applying \eqref{eq:smoothing} once more for $t\ge\tau$ gives
$\|v(t)\|_\infty\le C_\tau h_*^{-1}\|p(t-\tau)\|_2\le CD_0e^{-\lambda_1t}$.
For shorter times, enlarge the constant using \eqref{eq:kato}. This is precisely \eqref{eq:main1}.
\end{proof}

\section{Nonlinear spectral asymptotics and the proof of Theorem~\ref{thm:main2}}\label{sec:asymproof}
Throughout this section, fix a small-data solution covered by the first main
theorem. Unit-window $L^2(0,1;L^2(\Om))$ bounds for $\NN(\vct p)$ are used in
the spectral-semigroup, Bochner-integral, and energy arguments. The forcing
$g=-(\vct v\cdot\nabla)\varphi$ in the equation for the original-variable
remainder $Z$ instead has a pointwise-in-time $L^\infty(\Om)$ bound. This
permits a comparison estimate for its modulus. 

\subsection{An inhomogeneous comparison estimate}
\begin{lemma}[$L^\infty$ comparison with bounded drift]\label{lem:inhomogeneous}
Let $I=[s,s+\tau]$ with $0<\tau\le1$, and let $Z$ be the variational
solution of
\[
\begin{aligned}
Z_t-\mu\Delta Z+B(t,x)\cdot\nabla Z
   &=f-A(t,x)Z &&\text{in }\Omega,\\
Z|_{\partial\Omega}&=0,\qquad Z(s)=Z_s\in L^2(\Omega).
\end{aligned}
\]
Assume $\|B\|_\infty\le B_*$,
$\|A\|_{\infty,{\rm op}}\le M_*$,
and $f\in L^2(I;L^2)\cap L^1(I;L^\infty)$. Then, for each $t\in(s,s+\tau]$,
\begin{equation}\label{eq:inhomogeneous-main}
\|Z(t)\|_\infty\le C_{n,\mu}(t-s)^{-n/4}
 e^{(M_*+B_*^2/\mu)(t-s)}\|Z_s\|_2
 +\int_s^t e^{M_*(t-r)}\|f(r)\|_\infty dr.
\end{equation}
\end{lemma}
\begin{proof}
By linearity, decompose $Z=z+y$, where $z$ solves the homogeneous problem with initial value $Z_s$ and $y$ solves the inhomogeneous problem with zero initial value. Lemma~\ref{lem:moser} controls $z$. We establish the estimate for $y$ directly in the variational class, without requiring smooth approximation in the $L^1_tL^\infty_x$ topology.

Set $g(t)=\|f(t)\|_\infty\in L^1(I)$ and, for $\varepsilon>0$, define
\[
 \beta_\varepsilon(\xi)=\sqrt{|\xi|^2+\varepsilon^2}-\varepsilon,
 \qquad |D\beta_\varepsilon|\le1,\quad D^2\beta_\varepsilon\ge0.
\]
Lemma~\ref{lem:weak-kato} gives the convex-chain inequality in the variational class:
\begin{equation}\label{eq:direct-kato-regularized}
 (\partial_t-\mu\Delta+B\cdot\nabla)\beta_\varepsilon(y)
 \le M_*\beta_\varepsilon(y)+g(t)+M_*\varepsilon
 \quad\text{in distributions on }I\times\Om.
\end{equation}
Indeed, the diffusion contribution omitted on the right is nonpositive:
\[
 -\mu\sum_j D^2\beta_\varepsilon(y)[\partial_j y,\partial_j y]\le0.
\]
Moreover,
$|D\beta_\varepsilon(y)\cdot Ay|\le M_*|y|\le M_*(\beta_\varepsilon(y)+\varepsilon)$.
The trace of $\beta_\varepsilon(y)$ is zero since $y\in L^2(I;H_0^1)$ and $\beta_\varepsilon(0)=0$.

Let $c_\varepsilon$ solve
\[
 c_\varepsilon'=M_*c_\varepsilon+g+M_*\varepsilon,
 \qquad c_\varepsilon(s)=0,
\]
and put $Y_\varepsilon=(\beta_\varepsilon(y)-c_\varepsilon)_+$. Since
$0\le Y_\varepsilon\le\beta_\varepsilon(y)\le|y|$, the positive part belongs to
$L^2(I;H_0^1(\Om))\cap L^\infty(I;L^2(\Om))$; the spatially constant
$c_\varepsilon$ itself is \emph{not} required to belong to $L^2(\Om)$.
Testing the scalar inequality for $\beta_\varepsilon(y)-c_\varepsilon$ against
$Y_\varepsilon$, first on compact space-time truncations and then by density, gives
\[
 \frac12\frac{d}{dt}\|Y_\varepsilon(t)\|_2^2
 +\mu\|\nabla Y_\varepsilon(t)\|_2^2
 \le B_*\|Y_\varepsilon(t)\|_2\|\nabla Y_\varepsilon(t)\|_2
 +M_*\|Y_\varepsilon(t)\|_2^2.
\]
Using Young's inequality,
\begin{equation}\label{eq:comparison-positive-part}
 \frac12\frac{d}{dt}\|Y_\varepsilon\|_2^2
 +\frac\mu2\|\nabla Y_\varepsilon\|_2^2
 \le\left(M_*+\frac{B_*^2}{2\mu}\right)\|Y_\varepsilon\|_2^2.
\end{equation}
Here $Y_\varepsilon(s)=0$ in the strong $L^2$ trace. Steklov time averages and positive-part truncations justify the time chain rule used in \eqref{eq:comparison-positive-part}; the spatial cutoff terms vanish because $Y_\varepsilon\in L^2_tH_0^1\cap L^\infty_tL^2$ and $B$ is bounded. Gronwall's lemma gives $Y_\varepsilon=0$.

The inequality $\beta_\varepsilon(y(t))\le c_\varepsilon(t)$ holds at each fixed time because $y\in C(I;L^2)$ and the positive-part map is continuous in $L^2$ under continuous variation of its nonnegative scalar threshold. Passing to $\varepsilon\downarrow0$ therefore yields
\[
 \|y(t)\|_\infty\le c(t):=\int_s^t e^{M_*(t-r)}g(r)\,dr,
 \qquad s\le t\le s+\tau.
\]
Combining this estimate with the homogeneous smoothing bound proves \eqref{eq:inhomogeneous-main} without an assumption on $\divv B$; the bound holds at every fixed terminal time.
\end{proof}

\begin{proof}[Proof of Theorem~\ref{thm:main2}]
\emph{Step 1: unit-window forcing and the nonlinear amplitude.}
The energy absorption inequality \eqref{eq:half-energy}, the sharp decay estimate,
and the form-domain control \eqref{eq:form-coercive} from the proof of the first main theorem imply, for every real $t\ge0$,
\begin{equation}\label{eq:window}
\int_t^{t+1}\|p(s)\|_{H^1}^2ds\le CD_0^2e^{-2\lambda_1t},\qquad
\|\NN(p)\|_{L^2(t,t+1;L^2)}\le CD_0^2e^{-2\lambda_1t}.
\end{equation}
The second inequality follows from
$\|\NN(p)\|_2\le\|v\|_\infty(\|\nabla p\|_2+\|F\|_\infty\|p\|_2)$.
These bounds include the first window $[0,1]$ and do not assume smallness of $\nabla v_0$.
Consequently,
\[
\sum_{k\ge0}\int_k^{k+1}e^{\lambda_1s}\|\Pn\NN(p(s))\|_2ds
\le CD_0^2\sum_{k\ge0}e^{-\lambda_1k}<\infty.
\]
Thus the integral in \eqref{eq:amplitude} converges in $E_1$ as a Bochner integral, and
\begin{equation}\label{eq:A-bound}
\|\A(p_0)-\Pn p_0\|_2\le CD_0^2,\qquad \|\A(p_0)\|_2\le CD_0.
\end{equation}

\emph{Step 2: spectral decomposition and the three remainder regimes.}
Since $\Pn$ is a bounded spectral projection commuting with $e^{-tH}$, it may be moved inside
the Bochner integral in \eqref{eq:duhamel}, giving
\begin{equation}\label{eq:p-tail}
\Pn p(t)-e^{-\lambda_1t}A_\infty
=-e^{-\lambda_1t}\int_t^\infty e^{\lambda_1s}\Pn\NN(p(s))\,ds,
\end{equation}
as well as $\Qn p(t)=e^{-tH}\Qn p_0+\int_0^t e^{-(t-s)H}\Qn\NN(p(s))ds$.
We estimate the convolution over unit time windows without assuming pointwise continuity of the forcing.
If, in a Hilbert space $X$,
$\|g\|_{L^2(k,k+1;X)}\le Ge^{-\rho k}$, set $m=\lfloor t\rfloor$.
Applying Cauchy--Schwarz on each window yields
\[
\left\|\int_0^te^{-\Lambda(t-s)}g(s)ds\right\|_X
\le Ge^{-\Lambda t+\Lambda}\sum_{k=0}^{m-1}e^{(\Lambda-\rho)k}
+Ge^{-\rho m}.
\]
When $\Lambda<\rho$, $\Lambda=\rho$, or $\Lambda>\rho$, use respectively
a convergent geometric series, a sum of $m$ terms, or an increasing geometric series to obtain
\begin{equation}\label{eq:window-conv}
\left\|\int_0^te^{-\Lambda(t-s)}g(s)ds\right\|_X
\le C_{\rho,\Lambda}G
\begin{cases}
e^{-\Lambda t},&\Lambda<\rho,\\
(1+t)e^{-\rho t},&\Lambda=\rho,\\
e^{-\rho t},&\Lambda>\rho.
\end{cases}
\end{equation}
The argument for the last incomplete window also covers $0\le t<1$. For $0<\lambda<\rho$, the same window-by-window summation gives
\begin{equation}\label{eq:tail-conv}
\left\|e^{-\lambda t}\int_t^\infty e^{\lambda s}g(s)ds\right\|_X
\le C_{\lambda,\rho}G e^{-\rho t}.
\end{equation}
Apply \eqref{eq:tail-conv} to \eqref{eq:p-tail} with $\rho=2\lambda_1$ and $\lambda=\lambda_1$ to obtain
$\|\Pn p(t)-e^{-\lambda_1t}A_\infty\|_2\le CD_0^2e^{-2\lambda_1t}$.
On the orthogonal complement, $\|e^{-tH}\Qn\|_{2\to2}\le e^{-\Lambda t}$,
and thus \eqref{eq:window-conv} gives
\[
\|\Qn p(t)\|_2\le D_0e^{-\Lambda t}+CD_0^2\mathcal G_\Lambda(t).
\]
Combining the two estimates proves \eqref{eq:remainderL2}, also at $t=0$.
The factor $1+t$ agrees with the resonant convolution: for the scalar forcing $g(s)=e^{-\rho s}$, the convolution at
$\Lambda=\rho$ equals exactly $te^{-\rho t}$.

\emph{Step 3: eigenfunction regularity and the original-variable remainder.}
Take $\Psi\in E_1$. From
$-\mu\Delta\Psi=(\lambda_1\Id-\mtrx{V})\Psi$ and global Dirichlet elliptic regularity,
$\|\Psi\|_{H^2}\le C\|\Psi\|_2$.
To obtain a global gradient bound, compare the heat evolution with a bounded symmetric matrix potential to the scalar Dirichlet heat semigroup
using Kato's inequality: $|e^{-tH}g|\le e^{t\|\mtrx{V}\|_\infty}S_D(t)|g|$.
Thus the scalar heat-kernel $L^2\to L^\infty$ bound together with
$e^{-tH}\Psi=e^{-\lambda_1t}\Psi$ gives, at any fixed $t>0$,
$\|\Psi\|_\infty\le C_t\|\Psi\|_2$ in every dimension.
Moreover, for any $g\in L^2\cap L^\infty$, choose
$g_j\in C_c^\infty(\Om)$ such that $g_j\to g$ in $L^2$ and
$\sup_j\|g_j\|_\infty\le\|g\|_\infty+o(1)$.
Apply \eqref{eq:heat} and use convergence of $S_D(t)g_j\to S_D(t)g$ in the sense of distributional derivatives
to extend the short-time gradient estimate to this input class, without presupposing $g\in C_0$.
Here $\Psi,(\lambda_1\Id-\mtrx{V})\Psi\in L^2\cap L^\infty$.
For any fixed $0<t_0<1$, the eigenvalue equation has the Duhamel representation
\[
\Psi=S_D(t_0)\Psi+\int_0^{t_0}S_D(t_0-s)(\lambda_1\Id-\mtrx{V})\Psi\,ds.
\]
The short-time gradient estimate in \eqref{eq:heat} and
$\int_0^{t_0}(t_0-s)^{-1/2}ds=2\sqrt{t_0}$ yield
$\|\nabla\Psi\|_\infty\le C\|\Psi\|_2$.
Finite dimensionality then implies
\begin{equation}\label{eq:eigen-reg}
\begin{aligned}
&\|A_\infty\|_{H^2}+\|A_\infty\|_{W^{1,\infty}}\le CD_0,\\
&U:=h^{-1}A_\infty\in H^2\cap H_0^1\cap W^{1,\infty},\\
&\|U\|_{H^2\cap W^{1,\infty}}\le CD_0.
\end{aligned}
\end{equation}
The zero trace follows from $A_\infty\in H_0^1$. Since $U\in L^2\cap W^{1,\infty}$, it vanishes uniformly at infinity and thus has a representative in $C_0$.
For $\Psi\in\Dom(H)$, direct calculation in the local distributional sense, using $w=-2\mu F$, gives
\begin{equation}\label{eq:conjugacy}
h[-\mu\Delta+w\cdot\nabla+(\nabla w)](h^{-1}\Psi)
=-\mu\Delta\Psi+\mu(\divv F+|F|^2)\Psi-2\mu\nabla^2f\,\Psi=H\Psi.
\end{equation}
Therefore $\varphi(t)=e^{-\lambda_1t}U$ satisfies the linearized equation in the original variables and the homogeneous Dirichlet boundary condition.
Put $Z=v-\varphi$. Then
\begin{equation}\label{eq:Zeq}
Z_t-\mu\Delta Z+(w+v)\cdot\nabla Z+(\nabla w)Z=g,
\quad Z|_{\partial\Om}=0,
\quad g=-(v\cdot\nabla)\varphi.
\end{equation}
By the first main theorem, $h_*^{-1}<\infty$, and \eqref{eq:eigen-reg},
\begin{equation}\label{eq:mode-force}
\|g(t)\|_2+\|g(t)\|_\infty\le CD_0^2e^{-2\lambda_1t}
\quad\text{for every }t\ge0.
\end{equation}
In particular, no gradient smallness assumption is needed on the initial time window.
Moreover, $Z=h^{-1}[p-e^{-\lambda_1t}A_\infty]$ belongs to
$C_{\mathrm{loc}}L^2\cap L^2_{\mathrm{loc}}H_0^1$. Therefore \eqref{eq:Zeq}
can be treated by the nonhomogeneous comparison principle in the variational solution class described above.
By Lemma~\ref{lem:inhomogeneous}, for the same bounded drift and matrix term, the homogeneous contribution is controlled by short-time $L^2\to L^\infty$ smoothing, while the zero-initial-data forcing is controlled by modulus comparison. The forcing $g\in L_t^2L_x^2\cap L_t^1L_x^\infty$ is admissible, and the conclusion holds at every fixed time.
Thus, for fixed $0<\tau\le1$ and $t\ge\tau$,
\begin{equation}\label{eq:Zcomp}
\|Z(t)\|_\infty\le C_\tau\|Z(t-\tau)\|_2
+\int_{t-\tau}^t e^{M(t-r)}\|g(r)\|_\infty dr.
\end{equation}
Using $\|Z\|_2\le h_*^{-1}\|p-e^{-\lambda_1t}A_\infty\|_2$,
\eqref{eq:remainderL2} and \eqref{eq:mode-force},
and noting that shifting time by the fixed $\tau$ only modifies the constants preceding the exponentials, we obtain \eqref{eq:remainderLinfty}.
For $0\le t\le\tau$, use $\|v(t)\|_\infty+\|U\|_\infty\le CD_0$
and $e^{-\Lambda t}\ge e^{-\Lambda\tau}$ to complete the estimate. Thus the result includes time zero.

\emph{Step 4: unit-window gradient estimates.}
Write $R=p-e^{-\lambda_1t}A_\infty$. Then $R_t+HR=\NN(p)$, and
$R\in C_tL^2\cap L^2_{\rm loc}H_0^1$ and $R_t\in L^2_{\rm loc}H^{-1}$.
For any $\lambda_1<\eta<\min\{\Lambda,2\lambda_1\}$,
\eqref{eq:remainderL2} and \eqref{eq:gt} imply
$\|R(t)\|_2\le C_\eta D_0e^{-\eta t}$.
Test the weak equation against $R$ on each real interval $[t,t+1]$, estimating the bounded potential and forcing by Young's inequality. Then
\[
\int_t^{t+1}\|R(s)\|_{H^1}^2ds
\le C\left(\|R(t)\|_2^2+\int_t^{t+1}\|R(s)\|_2^2ds
+\int_t^{t+1}\|\NN(p(s))\|_2^2ds\right).
\]
Substituting \eqref{eq:window} and using $D_0\le1$, since $\eta<2\lambda_1$, yields \eqref{eq:remainderH1}.
The strong $L^2$ time trace and $L^2_tH^{-1}$ time derivative justify the Hilbert-triple energy identity at both endpoints of every window, including $t=0$. No initial gradient bound is required.

\emph{Step 5: the projected quadratic amplitude and its cubic remainder.}
Since $E_1$ is finite dimensional and \eqref{eq:eigen-reg} holds for each unit basis vector,
for $\Psi\in E_1$ and $a,b\in C_c^\infty$, integration by parts in only the terms containing $\nabla b$ gives the componentwise identity
\begin{equation}\label{eq:projected}
\begin{aligned}
2\langle\BB(a,b),\Psi\rangle
=\int_\Om h^{-1}\bigl[& (\divv a)(b\cdot\Psi)+(a\cdot\nabla\Psi)\cdot b\\
&-(b\cdot\nabla a)\cdot\Psi+(F\cdot b)(a\cdot\Psi)\bigr]dx.
\end{aligned}
\end{equation}
The term generated by $\nabla h^{-1}=-h^{-1}F$ cancels exactly with
the corresponding positive term in \eqref{eq:bilin}.
Consequently, the projected pairing with values in $E_1$ extends uniquely and continuously to $H_0^1\times L^2$ (and, by symmetry, to the reversed space), with
\begin{equation}\label{eq:Bbound}
\|\Pn\BB(a,b)\|_2\le C\|a\|_{H^1}\|b\|_2,\qquad
\|\Pn\BB(a,b)\|_2\le C\|a\|_2\|b\|_{H^1}.
\end{equation}

To make the constant explicit, choose an orthonormal basis
$\{\Psi_j\}_{j=1}^m$ of $E_1$. Since $\Psi_j\in W^{1,\infty}$ and
$h^{-1},F\in L^\infty$, each summand satisfies
$|\langle\BB(a,b),\Psi_j\rangle|\le C_j\|a\|_{H^1}\|b\|_2$.
Hence
\[
\|\Pn\BB(a,b)\|_2^2
 =\sum_{j=1}^m|\langle\BB(a,b),\Psi_j\rangle|^2
 \le\Bigl(\sum_{j=1}^m C_j^2\Bigr)\|a\|_{H^1}^2\|b\|_2^2.
\]
Because $C_c^\infty\times C_c^\infty$ is dense in $H_0^1\times L^2$,
the limit is independent of the approximating sequence; exchanging the arguments yields the second, symmetric bound.
The extension is defined by density at the $L^2$ regularity level and does not impose a boundary-trace condition.
This extension concerns only the finite-dimensional projected pairing; the full product $\BB(a,b)$ need not belong to $L^2$ in arbitrary dimensions.
In particular, the $L^2$ argument is not assigned a boundary trace.

For arbitrary $q\in L^2$, put $z(t)=S(t)q$ and choose
$\lambda_1/2<\beta<\lambda_1$. Applying the linear semigroup energy identity on $[\varepsilon,T]$
with the weight $e^{2\beta t}$ yields
\[
\frac12e^{2\beta T}\|z(T)\|_2^2+
\int_\varepsilon^T e^{2\beta t}\bigl(\QQ[z(t)]-\beta\|z(t)\|_2^2\bigr)dt
=\frac12e^{2\beta\varepsilon}\|S(\varepsilon)q\|_2^2.
\]
Since $\QQ[z]-\beta\|z\|_2^2\ge(1-\beta/\lambda_1)\QQ[z]\ge0$,
first let $\varepsilon\downarrow0$ and then $T\to\infty$. Strong continuity of the $L^2$ semigroup and the spectral lower bound give
\begin{equation}\label{eq:weighted-linear}
\sup_{t\ge0}e^{\beta t}\|S(t)q\|_2+
\|e^{\beta t}S(t)q\|_{L^2(0,\infty;H_0^1)}\le C_\beta\|q\|_2.
\end{equation}

The spectral measure $E_H(d\xi)$ also gives the exact form-domain dissipation identity from the initial time:
\[
\int_0^T\|H^{1/2}S(t)q\|_2^2dt
=\frac12\int_{[\lambda_1,\infty)}(1-e^{-2T\xi})\,d\|E_H(\xi)q\|_2^2
\le\tfrac12\|q\|_2^2.
\]
Thus, although the high-frequency linear gradient may be unbounded pointwise as $t\downarrow0$,
its time-integrated square is controlled by the $L^2$ initial data, independently of higher Sobolev norms of those data.

Accordingly, for arbitrary $q_1,q_2\in L^2$, \eqref{eq:Bbound}, Cauchy--Schwarz, and
$2\beta-\lambda_1>0$ prove absolute convergence of \eqref{eq:main2-quadratic} and the estimate
\begin{equation}\label{eq:A2bound}
\|\A_2(q_1,q_2)\|_2\le C_\beta\|q_1\|_2\|q_2\|_2.
\end{equation}
These estimates apply directly to $L^2$ initial data, with no $H_0^1$ hypothesis.

Finally, set $z(t)=S(t)p_0$ and $r(t)=p(t)-z(t)$.
The Duhamel formula \eqref{eq:duhamel}, \eqref{eq:window}, and the semigroup spectral bound give, for all $t\ge0$,
$\|r(t)\|_2\le CD_0^2e^{-\lambda_1t}$.
Test the difference equation $r_t+Hr=\NN(p)$ by energy methods on each integer window,
and then use \eqref{eq:window} and boundedness of the potential to obtain
\begin{equation}\label{eq:rwindow}
\begin{aligned}
&\|r\|_{L^\infty(k,k+1;L^2)}+\|r\|_{L^2(k,k+1;H^1)}
\le CD_0^2e^{-\lambda_1k},\\
&\|p\|_{L^\infty(k,k+1;L^2)}+\|p\|_{L^2(k,k+1;H^1)}\\
&\qquad+\|z\|_{L^\infty(k,k+1;L^2)}+\|z\|_{L^2(k,k+1;H^1)}
\le CD_0e^{-\lambda_1k}.
\end{aligned}
\end{equation}
The spectral theorem gives the window-energy estimate for the linear orbit, including the first window $k=0$.
By symmetry and bilinearity, $\NN(p)=\BB(p,p)$, and
$\BB(p,p)-\BB(z,z)=\BB(r,p)+\BB(z,r)$.
Combining \eqref{eq:Bbound} and \eqref{eq:rwindow} gives
\[
\int_k^{k+1}\|\Pn[\BB(p,p)-\BB(z,z)]\|_2dt
\le CD_0^3e^{-2\lambda_1k}.
\]
Multiplying by $e^{\lambda_1t}$ and summing over successive windows gives
\[
\begin{aligned}
&\left\|\A(p_0)-\Pn p_0-\A_2(p_0,p_0)\right\|_2\\
&\qquad\le CD_0^3\sum_{k=0}^\infty e^{-\lambda_1k}
\le\frac{C}{1-e^{-\lambda_1}}D_0^3,
\end{aligned}
\]
which is \eqref{eq:main2-Taylor}. The constants depend only on the fixed stationary state, the spectral lower
bound, the finite-dimensional projection, and the chosen small-data ball.
They are uniform over the mixed small-data ball, including data without
compact support or initial derivatives.
This proves Theorem~\ref{thm:main2}.
\end{proof}

\begin{corollary}[Sharpness of the lowest-mode decay rate]\label{cor:sharp-rate}
There are admissible smooth, compactly supported data for which
\begin{equation}\label{eq:sharp-limits}
\begin{split}
\lim_{t\to\infty}e^{\lambda_1t}\|h\vct v(t)\|_2&=\|\A(\vct p_0)\|_2>0,\\
\lim_{t\to\infty}e^{\lambda_1t}\|\vct v(t)\|_\infty
&=\|h^{-1}\A(\vct p_0)\|_\infty>0.
\end{split}
\end{equation}
Consequently, no exponent strictly larger than $\lambda_1$ can replace
$\lambda_1$ in \eqref{eq:main1}, uniformly over all admissible small data.
\end{corollary}
\begin{proof}
By density of smooth compactly supported fields in $L^2$ and $E_1\ne\{0\}$, there exists a field $\psi\in C_c^\infty(\Om;\R^n)$ with $\Pn\psi\ne0$. Let $p_0=\delta\psi$ and $v_0=\delta h^{-1}\psi$. The weight is smooth and positive, so $v_0$ remains smooth and compactly supported, and
$D_0=\delta(\|h^{-1}\psi\|_\infty+\|\psi\|_2)$ is sufficiently small for all $0<\delta<\delta_0$.
The amplitude bound \eqref{eq:A-bound} gives
\[
\A(\delta\psi)=\delta\Pn\psi+O_{L^2}(\delta^2),
\]
and therefore $\A(\delta\psi)\ne0$ for sufficiently small $\delta$.
By \eqref{eq:remainderL2} and \eqref{eq:remainderLinfty}, the corresponding
remainders tend to zero in $L^2$ and $L^\infty$ after multiplication by
$e^{\lambda_1t}$, since $\Lambda>\lambda_1$ and $2\lambda_1>\lambda_1$.
Continuity of the norms proves \eqref{eq:sharp-limits}. Since $h^{-1}>0$,
$h^{-1}\A(\delta\psi)$ is not the zero function, so both limits are strictly
positive. This does not assert a nonzero lowest-mode amplitude for every
nonzero initial perturbation.
\end{proof}

\begin{remark}[Sharp decay rate]
Choose $\psi\in C_c^\infty(\Omega;\mathbb R^n)$ with $\Pn\psi\ne0$.
Also choose any nonzero $\chi\in C_c^\infty(\Omega)$. Define
\[
\begin{aligned}
\vct p_{0,k}&=\delta\psi+c\delta h\chi\sin(kx_1)e_j,\\
\vct v_{0,k}&=h^{-1}\vct p_{0,k}.
\end{aligned}
\]
Choose $c>0$ so that $c\|h\chi\|_2\le\|\Pn\psi\|_2/4$. Then
$D_{0,k}\le C\delta$ uniformly in $k$ and
$\|\Pn\vct p_{0,k}\|_2\ge3\delta\|\Pn\psi\|_2/4$.
By~\eqref{eq:A-bound}, for fixed sufficiently small $\delta>0$,
$\|\A(\vct p_{0,k})\|_2\ge\delta\|\Pn\psi\|_2/2$, whereas
$\|\nabla\vct v_{0,k}\|_\infty\to\infty$. Thus the sharp exponent is attained
even among data with arbitrarily large initial gradients.
\end{remark}

\appendix
\section{Variational cutoff and time-trace details}\label{app:weighted-cutoffs}
The following localization result makes explicit the weak-limit argument used in Lemma~\ref{lem:localweighted}. It explains why bounded drift suffices and why differentiating the unbounded weight in the energy argument would be premature before establishing weighted integrability.

\begin{lemma}[Uniform weighted energy under truncation]\label{lem:truncation-principle}
Let $T<\infty$, and suppose $z\in C([0,T];L^2(\Om;\R^n))\cap L^2(0,T;H_0^1(\Om;\R^n))$ solves
\[
 z_t-\mu\Delta z+B(t,x)\cdot\nabla z+A(t,x)z=0
\]
variationally, with $B\in L^\infty((0,T)\times\Om)$ and $A\in L^\infty((0,T)\times\Om;\R^{n\times n})$. Let $a_N\in W^{1,\infty}(\Om)$ be nonnegative, bounded weights with
$|\nabla a_N|\le L a_N$ almost everywhere, where $L$ is independent of $N$. Then
\begin{equation}\label{eq:appendix-truncation}
 \sup_{0\le t\le T}\|a_Nz(t)\|_2^2
 +\mu\int_0^T\|a_N\nabla z(t)\|_2^2\,dt
 \le C_T\|a_Nz(0)\|_2^2,
\end{equation}
where $C_T$ does not depend on $N$. In particular, if $a_N=\min\{h,N\}$, then
\[
 \sup_{0\le t\le T}\|hz(t)\|_2^2+
 \int_0^T\|h\nabla z(t)\|_2^2\,dt
 \le C_T\|hz(0)\|_2^2
\]
whenever the initial weighted norm is finite.
\end{lemma}
\begin{proof}
Choose $\chi_R\in C_c^\infty(\overline\Om)$ with $\chi_R=1$ for $r\le R$,
$\chi_R=0$ for $r\ge2R$, and $|\nabla\chi_R|\le C/R$.
The test function $\chi_R^2 a_N^2z$ lies in $L^2(0,T;H_0^1)$ for each fixed $R,N$. To justify testing against a time-dependent $z$, first use Steklov time averages and then pass to the limit using the boundedness and Lipschitz regularity of $a_N$. Write
$Y_{N,R}=\|\chi_Ra_Nz\|_2^2$ and
$G_{N,R}=\|\chi_Ra_N\nabla z\|_2^2$. The diffusion--weight cross term is bounded by
\begin{align*}
 &2\mu\bigg|\sum_{i,j}\int\chi_R^2a_N(\partial_ja_N)z_i\partial_jz_i\,dx\bigg|
 \le\frac\mu8G_{N,R}+C\mu L^2Y_{N,R}.
\end{align*}
The cross term arising from the spatial cutoff satisfies
\begin{align*}
 &2\mu\bigg|\sum_{i,j}\int\chi_Ra_N^2(\partial_j\chi_R)z_i\partial_jz_i\,dx\bigg|\\
 &\qquad\le\frac\mu8G_{N,R}
 +\frac{C\mu}{R^2}\int_{R<r<2R}a_N^2|z|^2\,dx.
\end{align*}
The drift and matrix terms are estimated without integrating $B$ by parts:
\begin{align*}
 &\bigg|\int\chi_R^2a_N^2 z\cdot(B\cdot\nabla z)\,dx\bigg|
 \le\frac\mu8G_{N,R}+\frac{C\|B\|_\infty^2}{\mu}Y_{N,R},\\
 &\bigg|\int\chi_R^2a_N^2z\cdot Az\,dx\bigg|
 \le\|A\|_\infty Y_{N,R}.
\end{align*}
Consequently, for almost every time,
\begin{equation}\label{eq:appendix-cutoff-diff}
 \frac12Y_{N,R}'+\frac\mu2G_{N,R}
 \le C\Bigl(1+\mu L^2+\frac{\|B\|_\infty^2}{\mu}+\|A\|_\infty\Bigr)Y_{N,R}
 +\frac{C\mu}{R^2}\int_{R<r<2R}a_N^2|z|^2\,dx.
\end{equation}
For fixed $N$, $a_N$ is bounded and $z\in L^2(0,T;L^2)$, so the last term tends to zero in $L^1(0,T)$ as $R\to\infty$. Gronwall and Fatou therefore give \eqref{eq:appendix-truncation}, with the same constant for all $N$. For $a_N=\min\{h,N\}$, the bound
$|\nabla a_N|\le\|\nabla\log h\|_\infty a_N$ verifies the hypothesis. Applying monotone convergence separately to $|z|^2$ and $|\nabla z|^2$ completes the passage $N\to\infty$.
\end{proof}

\begin{remark}[Identification of the strong initial trace]\label{rem:trace-detail}
After the mild--variational identification, the preceding lemma gives
\[
\vct p=h\vct v\in L^\infty(0,T;L^2)\cap L^2(0,T;H_0^1).
\]
Since $\vct v$ is bounded,
$\NN(\vct p)\in L^2(0,T;L^2)$ and
\[
\vct p_t=-H\vct p+\NN(\vct p)\in L^2(0,T;H^{-1}).
\]
The Hilbert-triple theorem yields $\vct p\in C([0,T];L^2)$.
For each compactly supported test field $\vct\phi$,
\[
\langle \vct p(t),\vct\phi\rangle
=\langle \vct v(t),h\vct\phi\rangle.
\]
As $t\downarrow0$, the right-hand side converges to
$\langle h\vct v_0,\vct\phi\rangle$. The uniform weighted bound extends this
weak convergence to every $L^2$ test field. Hence the strong trace is
$\vct p(0)=h\vct v_0$; no initial weighted gradient bound is assumed.
\end{remark}

\section*{Acknowledgements}
Qingsong Zhao was supported by the National Natural Science Foundation of China
under Grant No. 12401281.

\section*{Statements and Declarations}
\subsection*{Author Contributions}
Qingsong Zhao: Conceptualization, Methodology, Formal analysis, Validation,
Writing -- original draft, Writing -- review and editing.

\subsection*{Competing Interests}
Qingsong Zhao declares that there are no known competing financial interests or
personal relationships that could have appeared to influence the work reported
in this paper.

\subsection*{Data Availability}
No datasets were generated or analyzed in this theoretical study.

\subsection*{Use of generative AI and AI-assisted technologies}
During the preparation of this manuscript, the author used OpenAI ChatGPT to
assist with language editing, literature comparison, the exploration and
refinement of mathematical arguments, and manuscript revision. The author
independently verified the mathematical arguments and references, reviewed and
revised the AI-assisted content, and takes full responsibility for the final
manuscript.


\begin{thebibliography}{99}

\bibitem{Burgers1948}
Burgers, J.M.: A mathematical model illustrating the theory of turbulence.
\emph{Adv. Appl. Mech.} \textbf{1}, 171--199 (1948).
\url{https://doi.org/10.1016/S0065-2156(08)70100-5}.

\bibitem{Hopf1950}
Hopf, E.: The partial differential equation $u_t+uu_x=\mu u_{xx}$.
\emph{Commun. Pure Appl. Math.} \textbf{3}, 201--230 (1950).
\url{https://doi.org/10.1002/cpa.3160030302}.

\bibitem{Cole1951}
Cole, J.D.: On a quasi-linear parabolic equation occurring in aerodynamics.
\emph{Quart. Appl. Math.} \textbf{9}, 225--236 (1951).
\url{https://doi.org/10.1090/qam/42889}.

\bibitem{Vaganan2012}
Vaganan, B.M.: Cole--Hopf transformations for higher dimensional Burgers
equations with variable coefficients. \emph{Stud. Appl. Math.} \textbf{129}(3), 300--308 (2012).
\url{https://doi.org/10.1111/j.1467-9590.2012.00551.x}.

\bibitem{Hashimoto2014}
Hashimoto, I.: Asymptotic behavior of radially symmetric solutions for the
Burgers equation in several space dimensions. \emph{Nonlinear Anal.} \textbf{100}, 43--58 (2014).
\url{https://doi.org/10.1016/j.na.2014.01.004}.

\bibitem{Hashimoto2016}
Hashimoto, I.: Behavior of solutions for radially symmetric solutions for Burgers
equation with a boundary corresponding to the rarefaction wave. \emph{Osaka J. Math.} \textbf{53}, 799--811 (2016).

\bibitem{HashimotoMatsumura2019}
Hashimoto, I., Matsumura, A.: Asymptotic behavior toward nonlinear waves for
radially symmetric solutions of the multi-dimensional Burgers equation.
\emph{J. Differential Equations} \textbf{266}, 2805--2829 (2019).
\url{https://doi.org/10.1016/j.jde.2018.08.045}.

\bibitem{YangZhaoZhao2021}
Yang, T., Zhao, H., Zhao, Q.: Asymptotics of radially symmetric solutions for the
exterior problem of multidimensional Burgers equation (in Chinese).
\emph{Sci. Sin. Math.} \textbf{51}, 1057--1072 (2021).
\url{https://doi.org/10.1360/SSM-2019-0204}.

\bibitem{ZhaoZhao2021}
Zhao, H., Zhao, Q.: Radially symmetric stationary wave for two-dimensional
Burgers equation. \emph{Discrete Contin. Dyn. Syst.} \textbf{41}(5), 2167--2185 (2021).
\url{https://doi.org/10.3934/dcds.2020357}.

\bibitem{Hashimoto2020}
Hashimoto, I.: Stability of the radially symmetric stationary wave of the Burgers
equation with multi-dimensional initial perturbations in exterior domain.
\emph{Math. Nachr.} \textbf{293}, 2348--2362 (2020).
\url{https://doi.org/10.1002/mana.201900233}.

\bibitem{BorchersMiyakawa1990}
Borchers, W., Miyakawa, T.: Algebraic $L^2$ decay for Navier--Stokes flows in
exterior domains. \emph{Acta Math.} \textbf{165}, 189--227 (1990).
\url{https://doi.org/10.1007/BF02391905}.

\bibitem{KozonoOgawa1994}
Kozono, H., Ogawa, T.: On stability of Navier--Stokes flows in exterior domains.
\emph{Arch. Ration. Mech. Anal.} \textbf{128}, 1--31 (1994).
\url{https://doi.org/10.1007/BF00380792}.

\bibitem{Henry1981}
Henry, D.: \emph{Geometric Theory of Semilinear Parabolic Equations}. \emph{Lecture Notes in
Mathematics}, vol. 840. Springer, Berlin (1981).
\url{https://doi.org/10.1007/BFb0089647}.

\bibitem{FoiasSaut1984}
Foias, C., Saut, J.-C.: Asymptotic behavior, as $t\to+\infty$, of solutions of
Navier--Stokes equations and nonlinear spectral manifolds. \emph{Indiana Univ. Math. J.} \textbf{33}, 459--477 (1984).
\url{https://doi.org/10.1512/iumj.1984.33.33025}.

\bibitem{FoiasSaut1987}
Foias, C., Saut, J.-C.: Linearization and normal form of the Navier--Stokes
equations with potential forces. \emph{Ann. Inst. H. Poincar\'e Anal. Non Lin\'eaire} \textbf{4}, 1--47 (1987).
\url{https://doi.org/10.1016/S0294-1449(16)30372-9}.

\bibitem{Schonbek1980}
Schonbek, M.E.: Decay of solution to parabolic conservation laws.
\emph{Commun. Partial Differ. Equ.} \textbf{5}, 449--473 (1980).
\url{https://doi.org/10.1080/0360530800882145}.

\bibitem{GhisiGobbinoHaraux2015}
Ghisi, M., Gobbino, M., Haraux, A.: A description of all possible decay rates for
solutions of some semilinear parabolic equations. \emph{J. Math. Pures Appl.} \textbf{103},
868--899 (2015).
\url{https://doi.org/10.1016/j.matpur.2014.09.005}.

\bibitem{HoangMartinez2017}
Hoang, L.T., Martinez, V.R.: Asymptotic expansion in Gevrey spaces for solutions
of Navier--Stokes equations. \emph{Asymptot. Anal.} \textbf{104}, 167--190 (2017).
\url{https://doi.org/10.3233/ASY-171429}.

\bibitem{Witten1982}
Witten, E.: Supersymmetry and Morse theory. \emph{J. Differential Geom.} \textbf{17}, 661--692
(1982).
\url{https://doi.org/10.4310/jdg/1214437492}.

\bibitem{HelfferNier2005}
Helffer, B., Nier, F.: \emph{Hypoelliptic Estimates and Spectral Theory for
Fokker--Planck Operators and Witten Laplacians}. \emph{Lecture Notes in Mathematics},
vol. 1862. Springer, Berlin (2005).
\url{https://doi.org/10.1007/b104762}.

\bibitem{LionsMagenes1972}
Lions, J.-L., Magenes, E.: \emph{Non-Homogeneous Boundary Value Problems and
Applications}, vol. I. Springer, Berlin (1972).
\url{https://doi.org/10.1007/978-3-642-65161-8}.

\bibitem{DenkHieberPruess2003}
Denk, R., Hieber, M., Pr\"uss, J.: $\mathcal R$-Boundedness, Fourier Multipliers
and Problems of Elliptic and Parabolic Type. \emph{Mem. Amer. Math. Soc.} \textbf{166}(788),
1--130 (2003).
\url{https://doi.org/10.1090/memo/0788}.

\bibitem{Lieberman1996}
Lieberman, G.M.: \emph{Second Order Parabolic Differential Equations}.
World Scientific, Singapore (1996).
\url{https://doi.org/10.1142/3302}.

\bibitem{Wood2007}
Wood, I.: Maximal $L^p$-regularity for the Laplacian on Lipschitz domains.
\emph{Math. Z.} \textbf{255}, 855--875 (2007).
\url{https://doi.org/10.1007/s00209-006-0055-6}.

\bibitem{FurutoIwabuchi2026}
Furuto, Y., Iwabuchi, T.: Higher-order derivative estimates for the heat equation
on a smooth domain. \emph{Commun. Pure Appl. Anal.} \textbf{36}, 55--72 (2026).
\url{https://doi.org/10.3934/cpaa.2026075}.

\bibitem{DLMF}
Olver, F.W.J., Olde Daalhuis, A.B., Lozier, D.W., Schneider, B.I., Boisvert, R.F.,
Clark, C.W., Miller, B.R., Saunders, B.V., Cohl, H.S., McClain, M.A. (eds.):
\emph{NIST Digital Library of Mathematical Functions},
Chapter 10: Bessel Functions, Version 1.2.8. National Institute of Standards and
Technology (2026). \url{https://dlmf.nist.gov/10} (accessed 25 November 2026).

\bibitem{Kato1995}
Kato, T.: \emph{Perturbation Theory for Linear Operators}, 2nd edn. Springer, Berlin
(1995).
\url{https://doi.org/10.1007/978-3-662-03930-2}.

\bibitem{Ouhabaz2005}
Ouhabaz, E.M.: \emph{Analysis of Heat Equations on Domains}. \emph{London Mathematical
Society Monographs Series}, vol. 31. Princeton University Press, Princeton (2005).

\end{thebibliography}
\end{document}